%% file: main.tex
\documentclass[mnsc,non-blindrev]{cls/informs3} 

\OneAndAHalfSpacedXI 

\usepackage{algorithm}

\EquationsNumberedThrough    

\TheoremsNumberedThrough     
\ECRepeatTheorems  %

\input{macros}

\bibpunct[, ]{(}{)}{,}{a}{}{,}%
\def\bibfont{\small}%
\allowdisplaybreaks
\begin{document}


\RUNAUTHOR{Ji, Pang and Xu }

\RUNTITLE{Online Optimization of Difference-of-Convex Compositions}

\TITLE{Online Optimization of Difference-of-Convex Compositions with Smooth Mappings}

\ARTICLEAUTHORS{%
\AUTHOR{Jingwei Ji}
\AFF{Management Science and Engineering, Stanford University, \EMAIL{jingwei.ji@stanford.edu}} 
\AUTHOR{Jong-Shi Pang}
\AFF{The Daniel J.\ Epstein Department of
Industrial and Systems Engineering, University of Southern California, \EMAIL{jongship@usc.edu}}
\AUTHOR{Renyuan Xu}
\AFF{Management Science and Engineering, Stanford University, \EMAIL{renyuanxu@stanford.edu}}
} 

\ABSTRACT{%
    We study online optimization for a broad class of structured non-convex non-smooth problems where each loss is a composition of a difference-of-convex function with a smooth mapping, and the feasible region is defined by constraint functions of the same kind. 
    We propose a time-smoothed proximal linear algorithm and a local-regret measure based on a proximal residual mapping. We show that this residual is a proper stationarity measure for the original problem: its fixed-point condition implies first-order stationarity. 
    Our analysis relies on a tangent-cone characterization for a feasible region described by composite difference-of-convex constraints, which is of independent interest and allows each update to be computed via a convex optimization oracle, despite the non-convexity of the problem. 
    We establish a local-regret bound and a bound on the total number of inner convex subproblems. 
    We also derive an error bound connecting the proximal residual to the distance to stationarity, providing a quantitative certificate of approximate stationarity. 
}%




\KEYWORDS{non-convex and non-smooth optimization, online optimization } 

\maketitle

\section{Introduction}

We study online learning for a broad class of structured non-convex non-smooth optimization problems, namely \underline{c}ompositions of \underline{d}ifference-of-\underline{c}onvex functions with \underline{s}mooth mappings (CDCS), with constraints defined by functions of the same kind. 
This introduces a state-of-the-art modeling framework
for structured non-convex non-smooth optimization into the online setting.

Though generic non-convex minimization is NP-hard, the optimization community has developed a rich body of research on structured non-convex non-smooth optimization. 
In particular, recent work has studied optimization problems in which the objective takes the form of a
CDCS function:
\begin{equation}
    x \mapsto(f-g) \circ \theta(x)=f(\theta(x))-g(\theta(x)),    
\end{equation}
where $f$ and $g$ are convex, and $\theta$ is smooth. 
Functions of this form are broad enough to cover a variety of applications, and yet are structured enough to admit efficient algorithms with convergence guarantees. 
They encompass a wide range of sparse and regularized learning models, where a smooth data-fitting term is coupled with a non-convex non-smooth penalty, including sparse reconstruction, variable selection, denoising, and inverse problems
\citep{wright2009sparse,zhang2010nearly,xu2021sparse,wang2018vector,le2015dc}. 
It also covers matrix and spectral optimization problems such as numerical-radius computation \citep{watson1996computing,he1997algorithm}. 
From a complementary regularization perspective, \cite{gao2024wasserstein} show that Wasserstein distributionally robust optimization induces variation regularization while accommodating possibly non-convex and non-smooth losses.
Recently, \cite{le2024minimizing} studies the offline minimization of this class of problems and develop a majorization-minimization algorithm with convergence guarantees to critical points.  

At the same time, online learning increasingly needs to handle non-convex objectives, because many sequential decision problems evolve in environments where convexity is not available. However, the classical external-regret benchmark is poorly suited to non-convex problems, since even the offline global minimization problem is generally intractable.
This motivates stationarity-based performance criteria, beginning with the local-regret framework of \cite{hazan2017efficient}. 
Subsequent work extended this line to dynamic and adaptive regret in non-stationary environments, to constrained and safe online non-convex optimization, and to settings with limited oracle feedback or composite regret; see, for example, \cite{aydore2018local, aydore2019dynamic, xu2024online, mulvaney2022dynamic, changdynamic, hallak2021regret, guan2023online, hallak2025sound}.

Despite the above advances, the existing online optimization literature does not cover the CDCS problem class studied in this work.
We consider an online setting where not only is the per-round loss of CDCS form, but also the feasible region is defined by CDCS constraints.  
This immediately creates two difficulties. First, a naive linearization of the problem need not preserve feasibility, so standard projected-gradient or proximal-gradient schemes are no longer directly applicable. Second, once the feasible region is non-convex, even formulating a computationally meaningful stationarity measure becomes nontrivial. Therefore, extending local-regret methods to this setting requires new geometric and algorithmic ideas.

\subsection{Contributions}

This work advances the frontiers of both of the aforementioned streams of literature.
We summarize our main contributions as follows. 

\begin{enumerate}
    \item We extend the local-regret framework of \citep{hazan2017efficient} to non-convex non-smooth optimization with both CDCS objectives and CDCS constraints. We introduce the regret measure in \eqref{eq:regret_def} and show that it is a proper measure of stationarity for the original problem, in the sense that the associated fixed-point condition implies stationarity, which is a
    necessary condition for local optimality; see Proposition~\ref{prop:fixed_point_to_stationarity_informal}. Furthermore, we propose a provably efficient algorithm and establish a local-regret bound for it in Theorem~\ref{thm:main}. 
    \item The tangent cone characterization in Theorem~\ref{thm:tangent_cone_characterization} for feasible regions described by CDCS constraints is of independent interest and may find use in contexts beyond the online learning problem studied here. Moreover, it is precisely this characterization that enables us to design an algorithm requiring only a convex optimization oracle, while still guaranteeing convergence to a stationary point.
    \item Based upon the above mentioned developments, we further establish a quantitative connection between the proposed proximal residual and the true first-order stationarity measure. In particular, in Proposition~\ref{prop:finite_error_bound}, we show that under suitable conditions, the distance from the updated point to the stationarity set can be upper bounded by the proximal residual. This provides a rigorous a posteriori justification for using the residual as a computationally tractable certificate of approximate stationarity.
\end{enumerate}

\subsection{Organization}

The rest of the paper is organized as follows.
Section~\ref{sec:setting} introduces the problem setting. In Section~\ref{subsec:pre_variational_analysis}, we review preliminaries from variational analysis. Section~\ref{subsec:cost_func_and_feasible_region} introduces the favorable CDCS class of loss functions and constraints. Section~\ref{subsec:local_regret} defines the time-smoothed objective, the proximal residual, and the associated local-regret measure, and Section~\ref{subsec:the_algorithm} presents the proposed online algorithm.
Section~\ref{sec:analysis} contains the main technical developments. Section~\ref{subsec:tangent_cone} establishes a tangent-cone characterization for feasible regions described by CDCS constraints, and Section~\ref{subsec:local_regret_guarantee} uses this result to justify the proposed residual as a proper stationarity measure and to establish the local-regret and iteration-complexity guarantees. Section~\ref{subsec:offline_implications} discusses the implications of the framework for offline non-convex non-smooth optimization, and Section~\ref{subsec:error_bound} derives an error bound relating the proximal residual to the distance to stationarity. The appendix contains the remaining proofs.

\subsection{Related work}

This work straddles two streams of literature: online learning and offline non-convex non-smooth optimization.

\textbf{Online learning.}
With the rapid development of online learning, there has been a surge of interest in the field to deal with non-convex problems. 
Overall, the literature in this area can be divided into two categories. Let $\ell_t$ be the loss function at time $t$. 
The first category focuses classical (external) regret bound 
\begin{equation} \label{eq:regret_def_global}
    \sum_{t=1}^T\ell_t\left(x_t\right)-\min _x \sum_{t=1}^T \ell_t(x)    
\end{equation}
and therefore, implicitly, global optimality. 
To make this goal attainable, the learner is often equipped with a powerful primitive, e.g., an offline optimization oracle \citep{suggala2020online, maillard2010online, xu2024online} or a sampling oracle \citep{krichene2015hedge}. 
Apart from the full information setting, the non-convexity has also been addressed in the stochastic bandit setting, see \cite{kleinberg2004nearly, kleinberg2008multi, bubeck2011x} for example.
In addition, works e.g. \cite{mulvaney2022dynamic, changdynamic} consider a dynamic regret notion, where the benchmark is a sequence of comparators that can change over time, and establish dynamic regret bounds for non-convex online learning.

Though this notion of external regret \eqref{eq:regret_def_global} is natural for online convex programming \citep{hazan2016introduction}, it is quite strong in nature for online non-convex programming. 
As for non-convex programming, even in the offline case, it is too ambitious to expect a polynomial time algorithm to achieve global optimality  \citep{murty1987some}. 
Undoubtedly, this makes these approaches less applicable in practice.

The second approach, pioneered by \cite{hazan2017efficient}, focuses on a more tractable notion of regret, called the \textit{local-regret}.
Roughly speaking, this criterion only requires the learner to output points that are approximately first-order stationary for the evolving objective.
Based upon this idea, extensions encompass, e.g., a class of non-convex non-smooth online programming \citep{hallak2021regret}, dynamic local-regret bounds \citep{aydore2019dynamic}.
\cite{guan2023online} consider a setting where the learner only has access to limited instantaneous oracle feedback. 
A more recent example includes \cite{hallak2025sound}, where the author develops a local-regret framework for online non-convex composite optimization based on a proximal-gradient stationarity measure. 
 
{ In contrast to this literature, our work develops a local-regret framework for online problems whose losses are CDCS functions and whose feasible region is itself defined by CDCS constraints. 
This requires a new proximal-linear residual and a new convex surrogate feasible set, because existing local-regret methods do not address the interaction between non-smooth DC-composite objectives and non-convex CDCS constraints.
} 

\textbf{Offline non-convex non-smooth optimization.} A vast body of variational-analysis research addresses generic non-convex, non-smooth objectives, yet a particularly tractable subclass is DC programming \citep{tao1997convex, hiriart1985generalized, moudafi2008difference, mainge2008convergence, tuy2009global}, where the objective is a difference of two convex functions. 
Recent advances in \cite{pang2017computing} shift attention to (DC-)constrained DC programming. They develop algorithms for computing d-stationary points over convex feasible sets and B-stationary points for a structured class of non-convex DC-constrained problems under a pointwise Slater condition. These guarantees are stronger than the DC-criticality generally obtained by classical DCA. In a related machine-learning application, \cite{tian2022computing} study the non-convex, non-smooth $\rho$-margin-loss SVM and develop a semi-proximal ADMM-based method that provably computes d-stationary points.
Building on this, \cite{lu2022penalty} further develop infeasible penalty algorithm that preserve B-stationarity under the same CQ while eliminating the need for an initial feasible point.
The class of convex composite optimization problems has been a focus of study due to its wide applicability \citep{burke1985descent, fukushima1981generalized, sagastizabal2013composite}. 
The strands of literature, composite programming and DC programming merge in the recent work of \cite{le2024minimizing}, who study the composition of DC functions with smooth mappings in the offline setting. 
Most recently, \cite{zhu2026nonconvex} study offline non-convex programs whose objective and functional inequality constraints have convex-composite structure. They analyze a smoothed prox-linear augmented-Lagrangian method under local regularity and establish finite-time guarantees in terms of a KKT residual.
{ In contrast to this offline literature, our setting imposes a smoothness condition on the second function in the DC decomposition, hence aiming for a stronger stationarity notion, and a finite-sample convergence guarantee. 
Moreover, we consider the more challenging CDCS constrained region, for which the tangent-cone characterization (cf. Theorem~\ref{thm:tangent_cone_characterization}) is the key technical ingredient that allows a convex proximal-linear surrogate to certify stationarity.  
This characterization may also be of independent interest beyond online optimization.
} 
Interested readers are referred to \cite{cui2021modern} for a comprehensive survey on the topic.


\section{Setting} \label{sec:setting}

In this section, we formalize the class of online optimization problems studied in this paper and introduce the main objects used in our analysis. We begin with several preliminaries from variational analysis, then present the favorable CDCS class of loss functions and constraints. Building on this structure, we define the time-smoothed objective, the associated proximal residual, and the local-regret measure, and finally introduce the proposed algorithm.

\subsection{Preliminaries from variational analysis}
\label{subsec:pre_variational_analysis}

\paragraph{Notation.}
Throughout the paper, the plus and minus signs between a point and a set denote the addition and subtraction of the point to each element of the set, respectively.
Namely, for any $x \in \mathbb{R}^n$ and $A \subseteq \mathbb{R}^n$, we have $x + A = \{x + a \mid a \in A\}$ and $x - A = \{x - a \mid a \in A\}$.
Throughout the paper, we use $\norm{\cdot}$ to denote the $\ell_2$ norm.
For a matrix $B$, we use $\sigma_{\min }^{+}(B)$ to denote the smallest nonzero singular value of $B$ (and 0 if rank of $B$ is zero). 
Given a function $F: \mathbb{R}^{n} \rightarrow \mathbb{R}^{m}$, we use $J F (x) $ to denote $\brackets{ \frac{\partial F_i}{\partial x_j} }_{i,j}$, the $m \times n$ Jacobian matrix. 
For a given set $X \subseteq \mathbb{R}^n$ and a vector $\bar{x} \in X$, the (Bouligand) tangent cone of $X$ at $\bar{x}$, denoted by $\mathcal{T}_{X}(\bar{x})$, is given by
$$
\mathcal{T}_{X}(\bar{x}) \defeq \left\{v \in \mathbb{R}^n \left\lvert\, \begin{array}{l}
\exists\left\{x^\nu\right\} \subset X \text { and }\left\{\tau_\nu\right\} \downarrow 0 \text { such that } \lim _{\nu \rightarrow \infty} x^\nu=\bar{x} \\
\text { and } \lim _{\nu \rightarrow \infty} \frac{x^\nu-\bar{x}}{\tau_\nu}=v
\end{array}\right.\right\} .
$$
Let $K$ be a cone in $\mathbb{R}^n$. The polar cone of $K$, denoted by $K^{*}$, is defined as
$
    K^{*} = \{ y \in \mathbb{R}^n \mid \langle y, x \rangle \leq 0, \ \forall x \in K \} .
$
Let $f: \mathbb{R}^n \to \mathbb{R} \cup \{+\infty\}$ be a lower semicontinuous convex function with 
$\mbox{Dom} f \defeq \left\{ x \in \mathbb{R}^n \left\lvert\, f(x) < \infty \right.\right\}$. 
At any point $x \in \operatorname{Dom} f$, the \emph{subdifferential} of $f$ is defined as
\begin{equation} \label{eq:subdifferential_def_convex}
    \partial f(x) = \{ v \in \mathbb{R}^n \mid f(y) \geq f(x) + v^{\top}(y - x), \ \forall y \in \operatorname{Dom} f \}.
\end{equation}
Let $\phi: \mathcal{O} \to \mathbb{R}$ be a \emph{locally Lipschitz} function defined on an open set $\mathcal{O} \subseteq \mathbb{R}^n$. The \emph{regular subdifferential} (or \textit{Fréchet subdifferential} ) of $\phi$ at a point $\bar{x} \in \mathcal{O}$, denoted $\widehat{\partial} \phi(\bar{x})$, is defined as
\begin{equation} \label{eq:regular_subdifferential_def}
    \widehat{\partial} \phi(\bar{x}) = \left\{ v \in \mathbb{R}^n \left\lvert\, \liminf_{x \to \bar{x}, \, x \neq \bar{x}} \frac{\phi(x) - \phi(\bar{x}) - v^{\top}(x - \bar{x})}{\|x - \bar{x}\|} \geq 0 \right. \right\}.
\end{equation}
The \emph{limiting subdifferential} of $\phi$ at $\bar{x}$, denoted $\partial \phi(\bar{x})$, is defined as
\begin{equation} \label{eq:limiting_subdifferential_def} \nonumber 
    \partial \phi(\bar{x}) = \left\{ v \in \mathbb{R}^n \mid \exists \{x^k\} \to \bar{x}, \ \{v^k\} \to v \text{ such that } v^k \in \widehat{\partial} \phi(x^k) \ \text{for all } k \right\}.
\end{equation}
When $\phi$ is convex, both the regular subdifferential and the limiting subdifferential coincide with the subdifferential defined in \eqref{eq:subdifferential_def_convex}.

Now consider a \emph{composite function} $g \circ \theta$, where $g: \mathbb{R}^m \to \mathbb{R} \cup \{+\infty\}$ is a lower semicontinuous convex function and $\theta: \mathbb{R}^n \to \mathbb{R}^m$ is a $C^1$ mapping. For any point $\bar{x} \in \operatorname{Dom}(g \circ \theta)$, we have the inclusion
\begin{equation} \label{eq:subdiff_chain_rule_regular}
    J\theta(\bar{x})^{\top} \partial g(\theta(\bar{x})) \defeq \left\{ J\theta(\bar{x})^{\top} y \mid y \in \partial g(\theta(\bar{x})) \right\} \subseteq \widehat{\partial}(g \circ \theta)(\bar{x}),
\end{equation}
as guaranteed by Theorem 10.6 of \cite{rockafellar2009variational}.

\subsection{The cost functions and the feasible region}
\label{subsec:cost_func_and_feasible_region}

In the literature \citep{cui2021modern}, a difference-of-convex (DC) function is called \textit{favorable} if it has a DC decomposition $f-g$, where $f$ and $g$ are both convex, and $g$ is differentiable.
In this work, we call a CDCS function \textit{favorable} if the following condition holds. 

\begin{definition} \label{def:regularity_condition}
    Given function triplet $f, g: \mathbb{R}^m \rightarrow \mathbb{R}$, $\theta: \mathbb{R}^{n} \rightarrow \mathbb{R}^m$ and set $\mathcal{A} \subseteq \mathbb{R}^n$, we say that $(f,g,\theta, \mathcal{A})$ satisfy the \textit{favorable condition} with parameters $(L_f, L_g, \beta_{\theta}, L_{\theta}, M)$ if the following conditions hold:
    \begin{itemize}[leftmargin=2em]
        \item $f$, $g$ are both proper convex;  
        \item $f$ is $L_f$-Lipschitz;
        \item $g$ is $L_g$-Lipschitz and differentiable;  
        \item $\theta$ is $\beta_{\theta}$-smooth and $ \norm{ J \theta } \leq L_{\theta}$ on $\operatorname{conv} (\mathcal{A})$; 
        \item ranges are bounded: $\abs{ f(\theta(x)) - g(\theta(x)) } \leq M $ for some constant $M$, for any $x \in \mathcal{A}$.
    \end{itemize}
\end{definition}

\vspace{0.5cm}

\paragraph{Per-round timeline and cost functions. }

We consider an online non-convex non-smooth optimization in an oblivious adversarial setting. 
During each iteration $t$, the learner is tasked to predict $x_t \in \mathcal{X} \subseteq \mathbb{R}^n$.
The feasible region $\mathcal{X}$ is fixed over time and is \textit{known} to the learner.  
We will specify the structure of $\mathcal{X}$ later in this section.
After taking the action $x_t$, the learner observes the entire cost function $\ell_t: \mathbb{R}^n \rightarrow \mathbb{R}$ (via access to a zeroth-order and a first-order oracle to $f_t, g_t, \theta_t$, to be defined shortly), and suffers a loss of $\ell_t(x_t)$.
The function $\ell_t$ is assumed to be a composition of a DC function and a smooth mapping, i.e.,
\begin{equation}
    \ell_t(x) = f_t( \theta_t ( x ) ) - g_t( \theta_t(x)), 
\end{equation}
where $f_t, g_t: \mathbb{R}^m \rightarrow \mathbb{R}$ and $\theta_t: \mathbb{R}^{n} \rightarrow \mathbb{R}^m$ satisfy the following blanket assumptions.
\begin{assumption} \label{assumption:regular_triple_obj}
    For each $t$, the function triplet $(f_t, g_t, \theta_t, \mathcal{X})$ satisfies the favorable condition in Definition~\ref{def:regularity_condition}, with uniform parameters $L_f, L_g, \beta_{\theta}, L_{\theta}, M$ over $t$.
\end{assumption}
In the case where $\beta_{\theta} = 0$, the problem reduces to the case of difference-of-convex optimization. 

In what follows, we provide several examples of CDCS functions that are commonly used in machine learning and signal processing.
First, we provide an example of a CDCS loss function that approximates the 0-1 loss used in classification problems.
\begin{example}[Smooth regression under approximated 0-1 loss] 
    \label{example:0-1_loss}
    Consider a binary classification problem, where the true label is $y \in \{-1, +1 \}$ and the feature vector is $a$. 
    The predicted label is generated by $\operatorname{sgn} \paren{ \theta(x, a) } \defeq \operatorname{sgn} \paren{ \theta_a(x) }$, where $\theta_a(x)$ is the predicted label with parameter $x$.
    The 0-1 loss is defined as $\ell(x) = \one{ \operatorname{sgn} \paren{ \theta_a(x) } \neq y} = \one{ \theta_a(x) \cdot y \leq 0 }$. 
    We may consider the following scaled ramp function to approximate the 0-1 loss:
    \begin{equation}
        \frac{1}{\delta} \paren{ \paren{ \delta-t }_+ - \paren{ -t }_+ }
        = \left\{
            \begin{array}{ll}
                1, & t < 0, \\
                1 - \frac{t}{\delta}, & 0 \leq t \leq \delta, \\
                0, & t > \delta,
            \end{array}
        \right.
    \end{equation}
    for some $\delta > 0$. 
    The function $r_{\delta}$ is a DC function, and $r_{\delta}(t) \rightarrow \one{t \leq 0}$ as $\delta \downarrow 0$ pointwise. 

    To get a smooth approximation, we can replace $(\cdot)_+$ by the softplus
    \begin{equation}
        \varphi_\tau(t):=\frac{1}{\tau} \log \left(1+e^{\tau t}\right) \quad \text{for} \quad \tau>0 . 
    \end{equation}
    Let $f(x) = \frac{1}{\delta} \varphi_\tau(\delta - x)$, $g(x) = \frac{1}{\delta} \varphi_\tau(-x)$, and $\theta_a(x) = \theta(x, a) \cdot y$. 
    Then the smooth approximation of the 0-1 loss can be expressed as a CDCS function:
    \begin{equation}
        \ell(x) = (f-g) \circ \theta_a (x) . 
    \end{equation}

    $\square$
\end{example}

\begin{example}[Smooth non-convex loss with convex regularization]
    Consider an online learning problem where, at round $t$, the learner incurs a loss of the form
    $
    \ell_t(x)=q_t(x)+r(x),
    $
    where $q_t: \mathbb{R}^n \rightarrow \mathbb{R}$ is a smooth, possibly non-convex loss function, and $r: \mathbb{R}^n \rightarrow \mathbb{R}$ is a convex regularizer. This class covers many standard formulations in machine learning and signal processing, where $q_t$ measures data-fitting error and $r$ promotes structure such as sparsity or robustness.

    We now show that \(\ell_t\) can be represented within the CDCS framework. 
    Define the smooth mapping
    $
        \theta_t(x)\defeq \bigl(x,\, q_t(x)\bigr)\in\mathbb{R}^{n+1}.
    $
    Next, define the outer functions
    $
        f_t(u,s)\defeq r(u)+s,
        g_t(u,s)\equiv 0,
        (u,s)\in\mathbb{R}^{n+1}.
    $
    Since \(r\) is convex and \(s\mapsto s\) is linear, the function \(f_t\) is jointly convex in \((u,s)\). 
    Moreover, \(g_t\equiv 0\) is convex and differentiable. 
    Composing with \(\theta_t\), we obtain
    \begin{equation}
        (f_t-g_t)\circ\theta_t(x)
        =
        f_t\bigl(x,q_t(x)\bigr)
        =
        r(x)+q_t(x)
        =
        \ell_t(x) , 
    \end{equation}
    which belongs to the CDCS class. 
$\square$
\end{example}

The following example illustrates that a finite maximum of convex or differentiable concave function composed with smooth mapping still yields a favorable CDCS function. 

 \begin{example}     \label{example:max_of_convex_or_concave}    
    Let $\varphi(x) = \max \braces{ h_i( \theta (x) ) }_{i=1}^{k}$, where each function $h_i$ is either a concave differentiable function, or a convex function (not necessarily differentiable). 
    The function $\theta$ is differentiable. 
    Then $\varphi$ is a favorable CDCS function; see also \cite{cui2021modern}, Section 6.1.3 for the underlying calculus. 
    $\square$ 
\end{example}

\vspace{0.5cm}

\paragraph{Feasible region. }
We aim to handle a broad non-convex constrained region, also characterized by a CDCS function. 
Namely, the feasible region is defined as
\begin{equation} \label{eq:feasible_region}
    \mathcal{X} \defeq \braces{ x \in \mathbb{R}^n \mid c(x) \defeq \phi( \zeta (x) ) - \psi( \zeta (x) ) \leq 0 } , 
\end{equation}
where $(\phi, \psi, \zeta, \mathbb{R}^n)$ satisfies the favorable condition with parameters $(L_{\phi}, L_{\psi}, \beta_{\zeta}, L_{\zeta}, \infty)$ in Definition~\ref{def:regularity_condition}.
The definition of $\mathcal{X}$  contains only one (scalar) constraint. 
We note that this is made for ease of exposition and entails no loss of generality.
This is because multiple constraints can be equivalently represented as a single constraint, as shown in the following lemma.  

\begin{lemma} \label{lemma:max_of_dc_is_dc}
    Assume for $i=1,2$:
    \begin{itemize}[leftmargin=2em]
        \item $f_i: \mathbb{R}^{m_i} \rightarrow \mathbb{R}$ is convex,
        \item $g_i: \mathbb{R}^{m_i} \rightarrow \mathbb{R}$ is convex and differentiable,
        \item $\theta_i: \mathbb{R}^n \rightarrow \mathbb{R}^{m_i}$ is differentiable.
    \end{itemize}
    Define $\phi_i(x)=\left(f_i-g_i\right) \circ \theta_i(x) = f_i\left(\theta_i(x)\right)-g_i\left(\theta_i(x)\right)$.
    Then
    $
    x \rightarrow \max \left\{\phi_1(x), \phi_2(x)\right\}
    $
    is a favorable CDCS function. In particular, there exist convex $f, g$ and a differentiable $\theta$ such that
    \begin{equation}
        \max \left\{\phi_1(x), \phi_2(x)\right\}=f(\theta(x))-g(\theta(x)),
    \end{equation}
    with $g$ differentiable.
\end{lemma}
The proof is provided in the Appendix.

Next, we provide an example demonstrating how a CDCS formulation can approximate an $\ell_0$ sparsity budget.

\begin{example}[Huber $\ell_0$ sparsity budget] 
    \label{example:sparsity_budget}

To this end, for given $\delta > 0$, we recall the Huber function \citep{gao1997waveshrink, gao1998wavelet} $H_\delta: \mathbb{R} \rightarrow \mathbb{R}^+$ defined as
\begin{equation}
    H_\delta(t)=
    \left\{
        \begin{array}{ll}
    \frac{t^2}{2 \delta}, & |t| \leq \delta, \\ [5pt]
    |t|-\frac{\delta}{2}, & |t|>\delta.
        \end{array}  
    \right.
\end{equation}
One can verify that $H_\delta$ is a convex and differentiable function. 

Given a vector $x \in \mathbb{R}^n$, consider the sparsity constraint
\begin{equation} \label{eq:sparsity_constraint_approximate_delta}
    \sum_{j=1}^{n} \one{ \abs{x_j} > \delta } \leq b 
\end{equation}
which enforces that at most $b$ entries of $x$ can exceed the threshold $\delta$ in magnitude.
To fulfill \eqref{eq:sparsity_constraint_approximate_delta}, it suffices to require that 
\begin{equation} \label{eq:sparsity_constraint_approximate_huber}
    \sum_{j=1}^{n} \paren{ \abs{x_j} - H_\delta(x_j) } \leq \frac{1}{2} \delta b . 
\end{equation}
To see this, we denote $I_{\delta} (x) = \braces{ j \in [n] : \abs{x_j} > \delta }$ and observe that 
\begin{equation} 
    \sum_{j \in I_{\delta} (x) } \frac{\delta}{2} + \sum_{j \notin I_{\delta} (x) } \underbrace{ \paren{ \abs{x_j} - \frac{x_j^2}{2\delta} } }_{ \geq 0  } 
    \geq \sum_{j \in I_{\delta} (x) } \frac{\delta}{2} 
    = \frac{1}{2} \delta \sum_{j=1}^{n} \one{ \abs{x_j} > \delta } . 
\end{equation}
To put \eqref{eq:sparsity_constraint_approximate_huber} in the form of \eqref{eq:feasible_region}, it suffices to set 
\begin{equation}
    \zeta(x)=x, \quad 
    \phi(z)=\sum_j\left|z_j\right|-\frac{1}{2} \delta b, \quad 
    \psi(z)=\sum_j H_\delta\left(z_j\right)  . 
\end{equation}
$\square$
\end{example}

\vspace{0.5cm}

\paragraph{Time-smoothed cost functions. }
To properly define the regret, as shown in, e.g., \cite{hazan2017efficient} and \cite{hallak2021regret}, it is necessary to consider a time-smoothed local version of the per-round cost function.
Otherwise, in an adversarial setting, it is impossible to achieve a sublinear regret.
The smoothing is over a sliding window of length $w \in \mathbb{N}$.
The time-smoothed cost function at time $t$ is defined as
\begin{equation}
    \Phi_t(x) \defeq \frac{1}{w} \sum_{s=t-w+1}^{t} \ell_s(x)
=  \frac{1}{w} \sum_{s=t-w+1}^{t} \brackets{ f_s( \theta_s ( x ) ) - g_s( \theta_s ( x ) )  } , 
\end{equation}
with the convention that $\ell_s(\cdot) = f_s = g_s = \theta_s = 0$ for $s \leq 0$.

\subsection{Local-regret}
\label{subsec:local_regret}

In view of \cite{hazan2017efficient}, instead of the regret notion studied in the majority of the online learning literature, namely to come up with predictions $x_t$'s that control the cumulative loss $\sum_{t=1}^{T} \Phi_t(x_t)$, we consider the \textit{local-regret}. 
This regret measures the non-stationarity of the smoothed cost functions $\Phi_t$'s, at the predictions $x_t$'s, as made clear below in what follows.

To define a proper measure of stationarity for the CDCS constrained optimization problem, we begin by defining the \textit{proximal linear} mapping at some point $x$.  
This mapping has to satisfy two requirements. First, it has to be practically computable, e.g., by solving a convex optimization problem. Second, it can provide a proper stationarity measure for the CDCS constrained optimization problem. 
To tackle the non-convexity of the feasible region,  we define the \textit{cushioned linearization set} at a reference point $\bar{x}$ 
by (i) linearizing the inner function $\zeta$ at 
$\bar{x}$ inside the convex function $\phi$,
(ii) linearizing the composite function $\psi \circ \zeta$ at $\bar{x}$, and (iii) adding and subtracting proximal terms to the two composite functions
$\phi \circ \zeta$ and $\psi \circ \zeta$, respectively.  This results in the following convex
set: 
\begin{eqnarray}  
    \tilde{\mathcal{L}}  (\bar{x}) &\defeq& \Bigg\{ x \in \mathbb{R}^n\,\, \Bigg|\,\,
    \phi \paren{ \zeta(\bar{x}) + J \zeta(\bar{x}) (x-\bar{x}) } 
    + \frac{1}{2} L_{\phi} \beta_{\zeta}  \norm{ x - \bar{x}}^2  \nonumber  \\
  &&\qquad \qquad \qquad   \leq \psi(\zeta(\bar{x}) ) +  \nabla \psi \paren{\zeta(\bar{x})}^{\top} J \zeta ( \bar{x} ) (x - \bar{x}) 
    - \frac{1}{2} L_{\psi} \beta_{\zeta}  \norm{ x - \bar{x}}^2  \label{eq:def_cushioned_linearization_set}
    \Bigg\},
\end{eqnarray}
which is nonempty as $ \bar{x} \in \tilde{\mathcal{L}}  (\bar{x})$.
Later in Theorem~\ref{thm:tangent_cone_characterization}, we will show this set captures the local geometry of the original non-convex feasible region $\mathcal{X}$. 
For the objective function, we apply the 
same operations to each term $f_s \circ \theta_s - g_s \circ \theta_s$, and obtain, with
$v_{t}^{g}(x) = \frac{1}{w}   \underset{s=t-w+1}{\overset{t}{\sum}}  J \theta_s (x)^{\top}  \nabla g_s ( \theta_s (x) )$, the \textit{proximal semi-linearized mapping}: for some scalar $\mu > 0$,
\begin{eqnarray} 
    S_{t}^{\mu} (x) \defeq S_{t}^{\mu} (x; v_{t}^{g}(x) )  &=&     \underset{z \in \tilde{\mathcal{L}}(x)  }{\argmin} 
    ~  \Big \{ \frac{1}{w} \sum_{s=t-w+1}^t ~ f_{s} \paren{ \theta_{s} (x) + J \theta_{s} (x) (z-x) } \nonumber \\ 
    & & - \inner{  v_{t}^{g}(x) , z - x} 
    + \frac{\mu}{2} \norm{z-x}^2 \Big \}  , \label{eq:prox_linear_constrained} 
\end{eqnarray}
The ``semi-linearization'' refers to the fact
that $f_s$ is not touched, only $g_s$ is linearized.
The \textit{proximal residual} mapping hence is defined as
\begin{equation} \label{eq:prox_gradient_mapping}
    G_{t}^{\mu} (x) \defeq \mu \cdot \brackets{ S_{t}^{\mu} (x ; v_{t}^{g}(x) ) - x } . 
\end{equation}
To clarify the motivation behind the above definition, particularly the use of the cushioned linearization set $\tilde{\mathcal{L}}(x)$ in the proximal semi-linearized mapping, we first present the following proposition informally.
The formal version, Proposition~\ref{prop:fixed_point_to_stationarity_formal}, and its proof will be given in Section~\ref{sec:analysis}.

\begin{proposition}[Informal] 
    \label{prop:fixed_point_to_stationarity_informal}
    Given $\bar{x} \in \mathcal{X}$. For  $\bar{v} = \frac{1}{w}  \underset{s=t-w+1}{\overset{t}{\sum}}  J \theta_s (\bar{x})^{\top}  \nabla g_s ( \theta_s (\bar{x}) )$, under suitable conditions, the fixed-point condition
    \begin{equation} \label{eq:fixed_point_informal}
        \bar{x} = S_{t}^{\mu}(\bar{x}; \bar{v})    
    \end{equation}
    implies that 
    \begin{equation} \label{eq:stationary_x_informal}
        0 \in \widehat{\partial} \Phi_t(\bar{x}) + \mathcal{T}_{\mathcal{X}} (\bar{x})^{*}  . 
    \end{equation}
\end{proposition}
Hence, the motivation for introducing the cushioned linearization set $\tilde{\mathcal{L}}(x)$ is that a fixed point of the associated proximal linear mapping \eqref{eq:fixed_point_informal} implies the stationary condition \eqref{eq:stationary_x_informal} for the original problem.
Indeed, \eqref{eq:stationary_x_informal} leads to Bouligand stationarity of $\Phi_t$ at $\bar{x}$ on $\mathcal{X}$, as shown in the following standard
lemma in variational analysis.
\begin{lemma} \label{lemma:zero_in_subdiff_implies_bouligand_stationarity}
    Let $X \subseteq \mathbb{R}^n$ be a closed set. 
    Given $\bar{x}$, suppose $\widehat{\partial} \phi (\bar{x})$ is nonempty.
    Let $\phi: \mathbb{R}^n \rightarrow \mathbb{R}$ be a Bouligand-differentiable function. 
    Then, we have 
    \begin{equation}
        0 \in \widehat{\partial} \phi (\bar{x}) + \mathcal{T}_{X} ( \bar{x} )^* 
        \Rightarrow \phi'( \bar{x}; d ) \geq 0 \text{ for all } d \in \mathcal{T}_{X} ( \bar{x} ) ,
    \end{equation}    
    which means that $\bar{x}$ is a Bouligand stationary point of $\phi$ on $X$.
\end{lemma}
The proof is standard, and is provided in Appendix for completeness.

Therefore, while the proximal linear mapping \eqref{eq:prox_linear_constrained} only requires a convex optimization problem to be solved, the fixed-point condition associated with it leads to a stationarity condition over a non-convex feasible region $\mathcal{X}$. 
Motivated by Proposition~\ref{prop:fixed_point_to_stationarity_informal}, it would seem natural to define the accumulative \textit{local-regret} to be 
\begin{equation} \label{eq:regret_def}
    \regret_T \defeq \sum_{t=1}^{T} \norm{ G_{t}^{\mu} (x_t) }^2 . 
\end{equation}
While $\operatorname{dist}\paren{0, \widehat{\partial}  \Phi_{t} (\cdot) + \mathcal{T}_{\mathcal{X}}(\cdot)^* } $ is also a natural candidate for a measure of stationarity, it is in general not computationally tractable.
Nevertheless, in Section~\ref{subsec:error_bound}, we will show that under suitable conditions, the proximal residual can be used to control the stationarity measure $\operatorname{dist}\paren{0, \widehat{\partial}  \Phi_{t} (\cdot) + \mathcal{T}_{\mathcal{X}}(\cdot)^* } $.
We refer  readers to Chapter 8 of \cite{cui2021modern} for a more thorough discussion.

\subsection{The algorithm}
\label{subsec:the_algorithm}

Now we are ready to present our algorithm.

\begin{algorithm}[H]
\caption{Time-smoothed online DC composite algorithm}
\label{alg:online_comp_dc_constrained}
\begin{algorithmic}[1]
\State \textbf{Input:} window size $w\ge1$, tolerance $\delta>0$, step size $\mu > 0$
\State Set $x_1 \in \mathcal{X}$ arbitrarily
\For{$t=1,\dots,T$}
    \State Predict $x_t$
    \State Observe the cost function $f_t$, $g_t$ and $\theta_t$; suffer the loss
    \State Initialize $k \gets 1$, $z_{t+1}^{k} \gets x_t$
    \While{$\norm{G_{t}^{\mu} ( z_{t+1}^{k} )} > \delta/w$}
        \State Compute
        \(
        v_{g}^{k} =
        \frac{1}{w}
        \sum_{s=t-w+1}^{t}
        J \theta_s (z_{t+1}^{k})^{\top}
        \nabla g_s \bigl( \theta_s (z_{t+1}^{k}) \bigr)
        \)
        \State Compute
        \begin{eqnarray}
                    z_{t+1}^{k+1} 
                &\leftarrow& 
                \underset{ z' \in \tilde{\mathcal{L}}( z_{t+1}^{k} ) }{\argmin} ~     
                    \Big \{ \frac{1}{w} \underset{s=t-w+1}{\overset{t}{\sum}}  
                    f_s \paren{  \theta_s (z_{t+1}^{k}) + J \theta_s (z_{t+1}^{k}) (z' - z_{t+1}^{k}) }  \nonumber \\
                & &
                    - \inner{ v_{g}^{k} , z' - z_{t+1}^{k} }
                    +  \frac{\mu}{2  } \norm{ z' - z_{t+1}^{k} }^2 \Big \}  \label{line:trial_d_opt_constrained}  
                \end{eqnarray}
        \State $k \gets k+1$
    \EndWhile
    \State $x_{t+1} \gets z_{t+1}^{k}$
\EndFor
\end{algorithmic}
\end{algorithm}


            



Algorithm~\ref{alg:online_comp_dc_constrained} implements the proximal linear framework introduced above in an online fashion. At each round \(t\), after observing the current loss function, the learner starts from the current prediction \(x_t\) and repeatedly applies the proximal linear mapping associated with the time-smoothed objective \(\Phi_t\). 
Each inner iteration solves a \textit{convex optimization problem} over the convex cushioned linearization set \(\tilde{\mathcal{L}}(z_{t+1}^k)\), and hence remains computationally tractable despite the non-convexity of the original problem. 
The inner loop terminates once the proximal residual becomes smaller than \(\delta/w\), which means that the output is approximately stationary with respect to the chosen residual measure. The resulting point is then used as the prediction for the next round.

\begin{remark}[Reduction to projected gradient \cite{hazan2017efficient}]
    Here we explain how our algorithm subsumes the projected-gradient method in \cite{hazan2017efficient}. 
    Let $\tilde{f}_t$ be a smooth convex function and $\mathcal{K}$ a convex set. 
    We set 
    \begin{equation}
        \ell_t(x)=f_t\left(\theta_t(x)\right), \quad f_t(y)=y, \quad g_t \equiv 0, \quad \theta_t(x)=\tilde{f}_t(x) . \nonumber 
    \end{equation}
    Then $f_t\left(\theta_t(x)+J \theta_t(x)(z-x)\right)=\tilde{f}_t(x)+\left\langle\nabla \tilde{f}_t(x), z-x\right\rangle$. 
    Hence, the proximal linear subproblem becomes 
    \begin{equation}
        S_{t}^{\mu} (x) 
        = \argmin_{z \in \mathcal{K}} ~ \braces{ 
        \left\langle\frac{1}{w} \sum_{s=t-w+1}^t \nabla \tilde{f}_s(x), z-x\right\rangle+\frac{\mu}{2}\|z-x\|^2 } 
        = \Pi_{ \mathcal{K} }  \paren{ x - \frac{1}{\mu} \frac{1}{w} \sum_{s=t-w+1}^{t}  \nabla \tilde{f}_s (x) } . \nonumber
    \end{equation}
    In this case, Algorithm~\ref{alg:online_comp_dc_constrained} reduces exactly to the time-smoothed online gradient descent algorithm proposed by \cite{hazan2017efficient}, generalized herein
    to a much broader context with a nonconvex nonsmooth composite objective function and a 
    nonconvex constraint set.
    
\end{remark}

\vspace{0.5cm}

As we will show in Section~\ref{sec:analysis}, this procedure enjoys both a local-regret guarantee and a finite bound on the total number of inner iterations.


\section{Analysis}
\label{sec:analysis}

This section contains the main theoretical developments of the paper.
In Section~\ref{subsec:tangent_cone}, we characterize the tangent cone of the feasible region under CDCS constraints and use this result to establish a formal version of Proposition~\ref{prop:fixed_point_to_stationarity_informal}. 
In Section~\ref{subsec:local_regret_guarantee}, we prove the main local-regret and iteration-complexity guarantees of the proposed algorithm. 
Section~\ref{subsec:offline_implications} discusses the implication of our framework for offline non-convex non-smooth optimization, while Section~\ref{subsec:computational_concerns} studies the existence of interior points for the convex surrogate subproblems. 
Finally, in Section~\ref{subsec:error_bound}, we derive an error bound relating the proximal residual to the distance to stationarity.

\subsection{Tangent cone characterization}
\label{subsec:tangent_cone}

In this section, we first characterize the tangent cone of the feasible region $\mathcal{X}$; the 
theory developed in this section can be of independent interest for the offline CDCS problems.
Then we provide a formal version of Proposition~\ref{prop:fixed_point_to_stationarity_informal}, which justifies the design of the proximal linear mapping and the associated local-regret. 

We recall from ~\eqref{eq:def_cushioned_linearization_set} that the cushioned linearization set at point $\bar{x}$ is defined to be 
\begin{eqnarray*}
    \tilde{\mathcal{L}}  (\bar{x}) &\defeq& \Bigg\{ x \in \mathbb{R}^n\,\, \Bigg|\,\,
    \phi \paren{ \zeta(\bar{x}) + J \zeta(\bar{x}) (x-\bar{x}) } 
    + \frac{1}{2} L_{\phi} \beta_{\zeta}  \norm{ x - \bar{x}}^2 \\
  &&\qquad \qquad \qquad   \leq \psi(\zeta(\bar{x}) ) +  \nabla \psi \paren{\zeta(\bar{x})}^{\top} J \zeta ( \bar{x} ) (x - \bar{x}) 
    - \frac{1}{2} L_{\psi} \beta_{\zeta}  \norm{ x - \bar{x}}^2 
    \Bigg\}. 
\end{eqnarray*}
This is a convex set as well. Moreover, it is nonempty if $\bar{x} \in \mathcal{X}$, since $ \bar{x} \in \tilde{\mathcal{L}}  (\bar{x})$ is a necessary and sufficient condition for $\bar{x} \in \mathcal{X}$.

The cushioned linearization set is designed to be convex, explicitly characterizable, and to preserve the local first-order geometry of the original non-convex set $\mathcal{X}$. 
There are two properties that are key to the success of its design. 
First, the two terms in the constraint function $c=\phi \circ \zeta-\psi \circ \zeta$ are linearized asymmetrically: for the convex term $\phi \circ \zeta$, we linearize the inner mapping $\zeta$, while for the concave term $-\psi \circ \zeta$, we linearize the whole composite expression. This asymmetric treatment is what ensures convexity of the resulting surrogate set, in the spirit of \cite{le2024minimizing}. 
Second, we add quadratic cushion terms to both sides of the linearized constraint. These terms compensate for the curvature induced by the smooth inner mapping and are crucial for ensuring that this set preserves the tangent-cone geometry of the original feasible region.
As we will show in Example~\ref{example:counter_example_without_cushion}, a naive design without the cushion terms may fail to preserve the tangent cone, even in a simple one-dimensional case.

As a natural candidate for the tangent cone of $\mathcal{X}$, we define
\begin{equation}
   C (\bar{x}) \defeq \braces{ d \in \mathbb{R}^n \mid 
   \phi^{\prime}(\zeta(\bar{x}) ; J \zeta(\bar{x}) d) \leq \nabla \psi(\zeta(\bar{x}))^{\top} J \zeta(\bar{x}) d }  . 
\end{equation}


\begin{definition}
    We say that the \textit{(CDCS) Slater constraint qualification (CQ) } holds at $\bar{x}$
    satisfying $\phi(\zeta(\bar{x})) = \psi(\zeta(\bar{x})$
    if there exists $d^0 \in \mathbb{R}^n$ such that 
    \begin{equation}
        \phi ' \paren{  \zeta ( \bar{x} ) ;  J \zeta ( \bar{x} ) d^0 }  < 
        \nabla \psi \paren{  \zeta ( \bar{x} ) }^{\top}  J \zeta ( \bar{x} ) d^0  . 
    \end{equation}
\end{definition}

We are now ready to present the main result of this section, which states that under the above Slater CQ, the tangent cone of the cushioned linearization set coincides with the tangent cone of the original feasible region, and both are given by $C(\bar{x})$.

\begin{theorem}[Tangent cone characterization] 
    \label{thm:tangent_cone_characterization}
        Consider some point $\bar{x} \in \mathcal{X}$ such that $c(\bar{x}) = 0$ and $\zeta(\bar{x}) \in \operatorname{int}( \operatorname{dom} \phi)$.
        We have the inclusion $\tilde{\mathcal{L}}  (\bar{x}) \subseteq \mathcal{X}$.
    Suppose there exists $d^0 $ satisfying 
    \begin{equation} \label{eq:constraint_qualification_d_0}
        \phi ' \paren{  \zeta ( \bar{x} ) ;  J \zeta ( \bar{x} ) d^0 }  < 
        \nabla \psi \paren{  \zeta ( \bar{x} ) }^{\top}  J \zeta ( \bar{x} ) d^0  . 
    \end{equation}

    Then
    \begin{equation}
        \mathcal{T}_{ \tilde{\mathcal{L}}  (\bar{x}) } ( \bar{x} )
    = \mathcal{T}_{ \mathcal{X} } (\bar{x}) 
    = C (\bar{x}) \ . 
    \end{equation}
\end{theorem}

\begin{myproof}
    
\begin{itemize}
    \item $ \mathcal{T}_{\mathcal{X}} ( \bar{x} )   \subseteq C(\bar{x}) $: 
    First, we claim that the function $c$ is Hadamard directionally differentiable at $\bar{x}$, and its directional derivative can be computed by the chain rule as follows
    \begin{equation}
        c^{\prime}(\bar{x} ; d)=\phi^{\prime}(\zeta(\bar{x}) ; J \zeta(\bar{x}) d)-\nabla \psi(\zeta(\bar{x}))^{\top} J \zeta(\bar{x}) d \nonumber . 
    \end{equation}
    To see this, we first consider the $\phi \circ \zeta$ part. 
    Because $\zeta(\bar{x}) \in \operatorname{int}( \operatorname{dom} \phi)$, the proper convex function $\phi$ is locally Lipschitz near $\zeta(\bar{x})$. This guarantees that its standard directional derivative coincides with its Hadamard directional derivative, making $\phi$ Hadamard directionally differentiable at $\zeta(\bar{x})$.  
    We note that differentiability implies Hadamard directional differentiability, and hence Bouligand differentiability.  
    Since both $\psi$ and $\zeta$ are differentiable, the composition $\psi \circ \zeta$ is Hadamard directionally differentiable, by virtue of Proposition 2.47 in \cite{bonnans2013perturbation}.

    Take any $ d \in \mathcal{T}_{\mathcal{X}} ( \bar{x} ) $. 
    By definition, there exists a sequence $t_k \downarrow 0$, $x^k \in \mathcal{X}$ with $x^k \rightarrow \bar{x}$ such that $d = \lim_{k \rightarrow \infty } ~ \frac{x^k - \bar{x}}{ t_k } $. 
    Define $d_k = \frac{x^k - \bar{x}}{ t_k }$. 
    Since $x^k \in \mathcal{X}$ and $c(\bar{x})=0$, we have for all $k$ that 
    $
        \frac{c\left(x^k\right)-c(\bar{x})}{t_k} \leq 0  . 
    $
    Now, since $c$ is Hadamard directionally differentiable at $\bar{x}$, we have 
    \begin{equation}
        \lim_{k \rightarrow \infty } \frac{c\left(\bar{x}+t_k d_k\right)-c(\bar{x})}{t_k} 
        = \phi^{\prime}(\zeta(\bar{x}) ; J \zeta(\bar{x}) d)-\nabla \psi(\zeta(\bar{x}))^{\top} J \zeta(\bar{x}) d  
        \nonumber . 
    \end{equation}
    Since $x^k=\bar{x}+t_k d_k$, every term on the left is nonpositive, so the limit is also nonpositive. Hence, we have $d \in C(\bar{x})$.

    \item $ \mathcal{T}_{ \tilde{\mathcal{L}}  (\bar{x}) } ( \bar{x} ) \subseteq \mathcal{T}_{ \mathcal{X} } (\bar{x}) $: 
    We note that $\tilde{\mathcal{L}}(\bar{x}) \subseteq \mathcal{X} $  for any $\bar{x} \in \mathcal{X}$ with $c(\bar{x})=0$. 
    Indeed, for any $x \in  \tilde{\mathcal{L}}  (\bar{x}) $, we can decompose 
    \begin{eqnarray}
        c(x) &=& \phi(\zeta(x))-\psi(\zeta(x)) \nonumber \\ 
        &=& \underbrace{\phi(\zeta(x)) - \phi \paren{ \zeta(\bar{x}) + J \zeta(\bar{x}) (x-\bar{x}) } }_{\text{Term (A)}} \nonumber \\ 
        & & + \underbrace{\psi(\zeta(\bar{x})) +  \nabla \psi(\zeta(\bar{x}))^\top J \zeta (\bar{x}) (x-\bar{x})  - \psi(\zeta(x)) }_{\text{Term (B)}} \nonumber \\ 
        & & + \underbrace{ \phi \paren{ \zeta(\bar{x}) + J \zeta(\bar{x}) (x-\bar{x}) }   
        - \psi(\zeta(\bar{x})) - \nabla \psi(\zeta(\bar{x}))^\top J \zeta (\bar{x}) (x-\bar{x}) }_{\text{Term (C)}} .
    \end{eqnarray}
    To bound Term (A), by Lipschitz continuity of $\phi$ and the smoothness of $\zeta$, we have
    \begin{equation} \label{eq:term_A}
       \text{Term (A)}
        \leq L_{\phi} \norm{ \zeta(x) -  \zeta(\bar{x}) - J \zeta(\bar{x}) (x-\bar{x}) }
        \leq \frac{1}{2} L_{\phi} \beta_{\zeta} \norm{ x - \bar{x} }^2 . 
    \end{equation}
    By convexity of $\psi$, we have
    \begin{equation}
        \psi (\zeta(x)) \geq \psi (\zeta(\bar{x})) + \nabla \psi(\zeta(\bar{x}))^\top  \paren{ \zeta(x) - \zeta(\bar{x}) }   . 
    \end{equation}
    Hence, Term (B) can be bounded as follows:
    \begin{eqnarray}
        \text{Term (B)} &\leq&  \nabla \psi(\zeta(\bar{x}))^\top J \zeta (\bar{x}) (x-\bar{x}) - \nabla \psi(\zeta(\bar{x}))^\top  \paren{ \zeta(x) - \zeta(\bar{x}) } \nonumber  \\ 
        &\leq&  \norm{ \nabla \psi(\zeta(\bar{x})) } \norm{ J \zeta (\bar{x}) (x-\bar{x}) -  \paren{ \zeta(x) - \zeta(\bar{x}) } } \label{eq:term_B_cauchy} \\ 
        &\leq& \frac{1}{2} L_{\psi} \beta_{\zeta} \norm{ x - \bar{x} }^2 \label{eq:term_B_smooth_lip} .
    \end{eqnarray}
    Inequality~\eqref{eq:term_B_cauchy} follows from the Cauchy-Schwarz inequality.  
    Inequality~\eqref{eq:term_B_smooth_lip} follows from the smoothness of $\zeta$ and the Lipschitz continuity of $\psi$.
    Finally, for $x \in \tilde{L}(\bar{x})$, Term (C) can be bounded by 
    \begin{equation} \label{eq:term_C}
        \text{Term (C)} \leq 
        - \frac{1}{2} L_{\phi} \beta_{\zeta} \norm{ x - \bar{x} }^2 
        - \frac{1}{2} L_{\psi} \beta_{\zeta} \norm{ x - \bar{x} }^2\  . 
    \end{equation}
    due to the definition of $\tilde{\mathcal{L}}(\bar{x})$. 
    Putting \eqref{eq:term_A}, \eqref{eq:term_B_smooth_lip} and \eqref{eq:term_C} together yields the fact that the condition $x \in \tilde{\mathcal{L}}(\bar{x})$ implies $x \in \mathcal{X}$.

    Therefore, the inclusion $\mathcal{T}_{ \tilde{\mathcal{L}}  (\bar{x}) } ( \bar{x} ) \subseteq \mathcal{T}_{ \mathcal{X} } (\bar{x}) $ follows from the standard fact that,  if set $A \subseteq B$ and $a \in A$, then the tangent cone $\mathcal{T}_A(a)$ is a subset of the tangent cone $\mathcal{T}_B(a)$.

    \item $ C(\bar{x})  \subseteq \mathcal{T}_{ \tilde{ \mathcal{L}} (\bar{x})  } (\bar{x}) $: 
    Take any $d \in C(\bar{x})$. 
    We aim to show that there exists a sequence of positive numbers $t_k \downarrow 0$ and a sequence of points $\tilde{x}^k \in  \tilde{ \mathcal{L}} (\bar{x}) $ such that $d = \lim_{k \rightarrow \infty} \frac{ \tilde{x}^k - \bar{x} }{ t_k } $.

    Let $\lambda_k = \frac{1}{k}$. 
    For each $k \geq 1$, we define 
    \begin{equation}
        u^{k} \defeq (1-\lambda_k) d + \lambda_k d^0 . 
    \end{equation}
    Clearly, we have $\lim_{k \rightarrow \infty} u^k = d$.
    Let $v^k \defeq J \zeta(\bar{x}) u^{k}$. 
    By the subadditivity and positive homogeneity of the directional derivative $\phi '(\zeta(\bar{x}); \cdot) $, we have 
    \begin{eqnarray}
        \phi '(\zeta(\bar{x}); v^{k} )  
        &\leq&   \phi '(\zeta(\bar{x}); (1-\lambda_{k} )  J \zeta(\bar{x}) d ) + \phi '(\zeta(\bar{x}); \lambda_{k}  J \zeta(\bar{x}) d^0 ) \nonumber \\ 
        &=&  (1-\lambda_k)  \phi '(\zeta(\bar{x});  J \zeta(\bar{x}) d ) +  \lambda_k \phi '(\zeta(\bar{x});  J \zeta(\bar{x}) d^0 )  \label{eq:phi_prime_leq} .
    \end{eqnarray}
     Recalling $ \phi^{\prime}(\zeta(\bar{x}) ; J \zeta(\bar{x}) d) \leq \nabla \psi(\zeta(\bar{x}))^{\top} J \zeta(\bar{x}) d $ since $d \in C(\bar{x})$ and the CQ condition \eqref{eq:constraint_qualification_d_0} for $d^0$, continuing from \eqref{eq:phi_prime_leq}, we have 
    \begin{eqnarray}
        \phi '(\zeta(\bar{x}); v^{k} ) &<&
        (1-\lambda_k) \nabla \psi(\zeta(\bar{x}))^{\top} J \zeta(\bar{x}) d   
        +  \lambda_k \nabla \psi \paren{  \zeta ( \bar{x} ) }^{\top}  J \zeta ( \bar{x} ) d^0    \nonumber \\ 
        &=&  \nabla \psi \paren{  \zeta ( \bar{x} ) }^{\top} v^{k}  .
    \end{eqnarray}
    Thus $\delta_k \defeq \nabla \psi \paren{  \zeta ( \bar{x} ) }^{\top} v^{k} - \phi ' \paren{ \zeta(\bar{x}); v^{k} }  > 0 $. 
    Moreover, we have 
    \begin{eqnarray}
        \delta_k &=& \nabla \psi \paren{  \zeta ( \bar{x} ) }^{\top} v^{k} - \phi ' \paren{ \zeta(\bar{x}); v^{k} }  \nonumber \\ 
        &\geq& (1-\lambda_k) \nabla \psi \paren{  \zeta ( \bar{x} ) }^{\top}  J \zeta ( \bar{x} ) d 
        +   \lambda_k \nabla \psi \paren{  \zeta ( \bar{x} ) }^{\top}  J \zeta ( \bar{x} ) d^0  \nonumber \\ 
        && - (1-\lambda_k)  \phi '(\zeta(\bar{x});  J \zeta(\bar{x}) d ) - \lambda_k \phi '(\zeta(\bar{x});  J \zeta(\bar{x}) d^0 ) \label{eq:delta_k_lower_bound_using_phi_prime_leq} \\
        &\geq& \lambda_k \paren{\nabla \psi \paren{  \zeta ( \bar{x} ) }^{\top}  J \zeta ( \bar{x} ) d^0 -  \phi^{\prime}(\zeta(\bar{x}) ; J \zeta(\bar{x}) d^0) } \label{eq:delta_k_lower_bound_d_in_C}   \\ 
        &=&  \lambda_k \delta_1 . 
        \label{eq:delta_k_lower_bound}
    \end{eqnarray}
    Inequality~\eqref{eq:delta_k_lower_bound_using_phi_prime_leq} follows from \eqref{eq:phi_prime_leq}. 
    Inequality~\eqref{eq:delta_k_lower_bound_d_in_C} follows from the fact that $d \in C(\bar{x})$.

    The function $\phi$ is convex and hence $\phi ' ( \zeta(\bar{x}); v^k)$ is well defined. Hence, 
    for each $k \geq 1$, there exists some $\tau_k$ such that 
    \begin{equation}
        \frac{ \phi( \zeta(\bar{x}) + s_k v^{k} ) - \phi(\zeta(\bar{x}) ) }{s_k} 
        < \phi ' \paren{ \zeta(\bar{x}); v^{k} } + \frac{1}{2} \lambda_k \delta_1 . 
    \end{equation}
   for all $s_k \in (0, \tau_k]$.
   Hence, multiplying both sides by $s_k$, we have
    \begin{eqnarray}
        \phi( \zeta(\bar{x}) + s_k v^{k} )  &<&  \phi(\zeta(\bar{x}) ) + s_k \phi ' \paren{ \zeta(\bar{x}); v^{k} } +  \frac{1}{2} s_k \lambda_k \delta_1  \nonumber \\ 
        &\leq&  \phi(\zeta(\bar{x}) ) + s_k ( \nabla \psi \paren{  \zeta ( \bar{x} ) }^{\top} v^{k} -    \lambda_k \delta_1 ) + \frac{1}{2} s_k \lambda_k \delta_1  \nonumber \\ 
        &=&  \phi(\zeta(\bar{x}) ) + s_k  \nabla \psi \paren{  \zeta ( \bar{x} ) }^{\top} v^{k}  -  \frac{1}{2} s_k \lambda_k \delta_1  ,  \label{eq:phi_zeta_bar_s_v_k}
    \end{eqnarray}
    where the second inequality follows from \eqref{eq:delta_k_lower_bound}. 
    Now, we set 
    \begin{equation}
        t_k = \min \braces{ \tau_k , ~ 
        \frac{ \lambda_k \delta_1 }{   (L_{\phi} + L_{\psi} ) \beta_{\zeta} \norm{u^k}^2 }, \frac{1}{k}  } . 
    \end{equation}
   We have 
   \begin{eqnarray}
       \phi(\zeta(\bar{x}) + t_k v^{k} )  &<&  \phi(\zeta(\bar{x}) ) + t_k  \nabla \psi \paren{  \zeta ( \bar{x} ) }^{\top} v^{k}  -  \frac{1}{2} t_k \lambda_k \delta_1  \label{eq:phi_zeta_bar_t_v_k_1} \\ 
        &\leq&   \psi(\zeta(\bar{x}) ) + t_k  \nabla \psi \paren{  \zeta ( \bar{x} ) }^{\top} v^{ k}  -  \frac{1}{2} (L_{\phi} + L_{\psi} ) \beta_{\zeta}  \norm{ t_k u^k}^2 \label{eq:phi_zeta_bar_t_v_k_2}  , 
   \end{eqnarray}
    which implies $\bar{x} + t_k u^k \in \tilde{\mathcal{L}}(\bar{x})$. 
   Inequality~\eqref{eq:phi_zeta_bar_t_v_k_1} follows from setting $s_k$ to be $t_k$ in \eqref{eq:phi_zeta_bar_s_v_k}. 
   Inequality~\eqref{eq:phi_zeta_bar_t_v_k_2} follows from the choice of $t_k$. 
   To finish the proof, it suffices to realize that we can take $\tilde{x}^k = \bar{x} + t_k u^k$.
\end{itemize}

\end{myproof}

If $\bar{x} \in \mathcal{X}$ satisfies $c(\bar{x}) < 0$, then $\bar{x}$ is in the interior of $\mathcal{X}$, and hence $\mathcal{T}_{\mathcal{X}} (\bar{x}) = \mathcal{T}_{\tilde{L}(\bar{x})}(\bar{x})  = \mathbb{R}^n$.

\begin{corollary} \label{corollary:first_order_tangent_cone}
    Let $ y \in \partial \mathcal{X}$, and assume that the CDCS Slater CQ is satisfied at $y$. Then for every $q \in \partial \phi(\zeta(y))$, we have 
    \begin{equation}
        J \zeta(y)^{\top}(q-\nabla \psi(\zeta(y))) \in \mathcal{T}_{\mathcal{X}}(y)^* . 
    \end{equation}
\end{corollary}

\begin{myproof}
    Take any $h \in \mathcal{T}_{\mathcal{X}}(y)$. By Theorem~\ref{thm:tangent_cone_characterization}, we have $h \in C(y)$, which implies $\phi '(\zeta(y); J \zeta(y) h) \leq \nabla \psi(\zeta(y))^{\top} J \zeta(y) h$. Since $q \in \partial \phi(\zeta(y))$, we have $\langle q, J \zeta(y) h\rangle \leq \phi^{\prime}(\zeta(y) ; J \zeta(y) h)$. Hence, we have $\left\langle J \zeta(y)^{\top}(q-\nabla \psi(\zeta(y))), h\right\rangle \leq 0$, which implies $J \zeta(y)^{\top}(q-\nabla \psi(\zeta(y))) \in \mathcal{T}_{\mathcal{X}}(y)^*$. 
\end{myproof}

\vspace{0.5cm}

In seeking the right way to characterize the tangent cone of the feasible set $\mathcal{X}$, one may consider the following linearization set
\begin{equation}
    \mathcal{L} (\bar{x}) \defeq \braces{ x \in \mathbb{R}^n \mid 
    \phi \paren{ \zeta(\bar{x}) + J \zeta(\bar{x}) (x-\bar{x}) } 
    \leq \psi(\zeta(\bar{x}) ) +  \nabla \psi \paren{\zeta(\bar{x})}^{\top} J \zeta ( \bar{x} ) (x - \bar{x}) } 
\end{equation}
which is clearly a convex set. 
This set is a natural candidate for DC problems, as considered in e.g., \cite{pang2017computing}.
However, in our setting the composition with the inner $\zeta$ makes the problem more challenging. In this case, the linearization set $\mathcal{L} (\bar{x})$ may not be a subset of $\mathcal{X}$, as shown in the example below.


\begin{example}[$\mathcal{L}(\bar{x})$ may not be a subset of $\mathcal{X}$] \label{example:counter_example_without_cushion}
    Consider the following specification on the real line: 
    $ \phi(x) = \log \paren{1+e^x} , \psi(x) = \log \paren{1+e^{-x}} , \zeta(x) = \frac{x^2-1}{1+x^2}$. 
    It is easy to verify that the regularity assumptions are satisfied. 
    Since $c(x) = \phi(\zeta(x)) - \psi(\zeta(x)) = \log \paren{  \frac{ 1+e^{\zeta(x)} }{  1+e^{-\zeta(x)} } }= \zeta(x) =  \frac{x^2-1}{1+x^2} $, the feasible set is $\mathcal{X} = \braces{ x \in \mathbb{R} \mid  \frac{x^2-1}{1+x^2} \leq 0 } = [-1, 1]$.
    Let the reference point be $\bar{x} = 1$, which is on the boundary of $\mathcal{X}$. 
    At this point,  $\zeta(1)=0, J \zeta(1)=\zeta^{\prime}(1)=1 $. So the naive linearization set is 
    \begin{equation}
        \mathcal{L}(1)=\{x: \phi(\zeta(1)+J \zeta(1)(x-1)) \leq \psi(\zeta(1))+\nabla \psi(\zeta(1)) J \zeta(1)(x-1)\}
        = \left\{x: \log \left(1+e^{x-1}\right) \leq \log 2-\frac{x-1}{2}\right\}
        \nonumber . 
    \end{equation}
    Let $h(s)=\log \left(1+e^s\right)+\frac{s}{2}-\log 2$. Then $h(0)=0$ and  $h^{\prime}(s)=\frac{e^s}{1+e^s}+\frac{1}{2}>0$ for all $s$. 
    Then we can conclude that $\mathcal{L}(1) = (-\infty, 1]$, implying that $\mathcal{L}(1) \nsubseteq \mathcal{X}$. 
    Moreover, the Slater condition translates to the existence of some $d^0$ such that $\phi^{\prime}(0) d^0<\psi^{\prime}(0) d^0$, which is clearly satisfied.
    $\square$
\end{example}

\vspace{1cm}

Now, we are ready to state and prove the formal version of Proposition~\ref{prop:fixed_point_to_stationarity_informal}. 

\begin{proposition}[Formal version of Proposition~\ref{prop:fixed_point_to_stationarity_informal}] \label{prop:fixed_point_to_stationarity_formal}
    Given $\bar{x} \in \mathcal{X}$. 
    Suppose either 
    \begin{enumerate}
        \item $c(\bar{x}) < 0$, or 
        \item $c(\bar{x}) = 0$ and $\zeta(\bar{x}) \in \operatorname{int} ( \operatorname{dom} \phi )$. Suppose the Slater's condition holds at $\bar{x}$, i.e., there exists $d^0$ satisfying 
        \begin{equation} 
            \phi ' \paren{  \zeta ( \bar{x} ) ;  J \zeta ( \bar{x} ) d^0 }  < 
            \nabla \psi \paren{  \zeta ( \bar{x} ) }^{\top}  J \zeta ( \bar{x} ) d^0  . 
        \end{equation}
    \end{enumerate}
    For $v = \frac{1}{w}  \underset{s=t-w+1}{\overset{t}{\sum}}  J \theta_s (\bar{x})^{\top}  \nabla g_s ( \theta_s (\bar{x}) )$,  the fixed-point condition
    \begin{equation}
        \bar{x} = S_{t}^{\mu}(\bar{x}; v)    
    \end{equation}
    implies that 
    \begin{equation} 
        0 \in \widehat{\partial} \Phi_t(\bar{x}) + \mathcal{T}_{\mathcal{X}} (\bar{x})^{*}  . 
    \end{equation}
\end{proposition}

\begin{myproof}
    
For ease of exposition, we omit the time index $t$.
Without loss of generality, we consider the special case $w=1$. 
Given $x$, we write 
$$
S(x; v) = \argmin_{z \in \tilde{\mathcal{L}}(x)} ~ \braces{ f \paren{ \theta(x) + J \theta(x) (z-x) } - \inner{  v , z-x } + \frac{\mu}{2} \norm{ z-x }^2 } 
$$
for some positive $\mu$.
The objective in the curly brackets is strongly convex in $z$, so the minimizer $S(x; v)$ is unique and satisfies the first-order condition
$$
0 \in J \theta(x)^{\top} \partial f \paren{ \theta(x) + J \theta(x) (S(x; v)-x) } - v + \mu ( S(x; v) - x) + \mathcal{T}_{\tilde{\mathcal{L}}(x)} ( S(x; v) )^{*}  . 
$$
For $v = J \theta( \bar{x} )^{\top} \nabla g(\theta(\bar{x}))$, the fixed-point condition $ \bar{x} = S(\bar{x}; v) $ leads to 
\begin{equation} \nonumber 
    0 \in J \theta(\bar{x})^{\top} \partial f \paren{ \theta(\bar{x}) } - J \theta( \bar{x} )^{\top} \nabla g(\theta(\bar{x})) 
+ \mathcal{T}_{\tilde{\mathcal{L}}(\bar{x})} ( \bar{x} )^{*} . 
\end{equation}
For the case when $c(\bar{x}) < 0$, we have $\mathcal{T}_{\tilde{\mathcal{L}}(\bar{x})} ( \bar{x} )^{*} = \{0\}$. 
For the case when $c(\bar{x}) = 0$, then by Theorem~\ref{thm:tangent_cone_characterization} and since the Slater's condition holds, we can translate the above condition to 
\begin{eqnarray}
    0 &\in& J \theta(\bar{x})^{\top} \partial f \paren{ \theta(\bar{x}) } - J \theta( \bar{x} )^{\top} \nabla g(\theta(\bar{x})) 
+ \mathcal{T}_{\mathcal{X} } ( \bar{x} )^{*} \nonumber \\ 
    &\subseteq& \widehat{\partial} ( f \circ \theta ) (\bar{x}) - \widehat{\partial} ( g \circ \theta ) (\bar{x}) 
    +  \mathcal{T}_{\mathcal{X} } ( \bar{x} )^{*} \label{eq:apply_comp_subgradient_twice} \\ 
    &\subseteq& \widehat{\partial} \paren{ f \circ \theta - g \circ \theta } (\bar{x}) 
    +  \mathcal{T}_{\mathcal{X} } ( \bar{x} )^{*} .  \nonumber 
\end{eqnarray}
The inclusion~\eqref{eq:apply_comp_subgradient_twice} follows by applying the fact \eqref{eq:subdiff_chain_rule_regular} twice.
The last inclusion follows from Proposition~1.107 in \cite{mordukhovich2009variational}, as $ g \circ \theta $ is differentiable. 

The proof is now complete. 
\end{myproof}

\subsection{Local-regret guarantee }
\label{subsec:local_regret_guarantee}

We are now ready to state the main performance guarantee for Algorithm~\ref{alg:online_comp_dc_constrained}.
The theorem below shows that the proposed method attains a local-regret bound of order $\frac{T}{w^2}$, which is tight in time-smoothed local-regret analysis \citep{hazan2017efficient}.
At the same time, it controls the total amount of optimization-oracle calls required by the inner proximal linear loop.

\begin{theorem}  \label{thm:main}
    Let $\ell_1, \cdots, \ell_T$ be the sequence of cost functions presented to Algorithm~\ref{alg:online_comp_dc_constrained}.
    Suppose Assumption~\ref{assumption:regular_triple_obj} holds, then
    \begin{enumerate}[label=(\roman*)]
        \item Given $\mu >  (L_f + L_g) \beta_{\theta} $, the local-regret incurred is bounded by 
        \begin{equation} \label{eq:regret_bound}
        \regret_T \lesssim T \cdot \paren{ \frac{\delta + (L_f+L_g) L_{\theta} }{w} }^2 . 
        \end{equation}
        \item  Let $\tau_t$ be the number of iterations of the while loop executed  at time $t$.  
        The total number of inner optimization iterations is bounded by 
        \begin{equation}
            \sum_{t=1}^{T} \tau_t \lesssim 
             \frac{ M \mu^2 }{\mu-\left(L_f+L_g\right) \beta_\theta} \frac{w^2}{\delta^2} \left(1+\frac{T}{w}\right)  . 
        \end{equation}
    \end{enumerate}

\end{theorem}

\begin{myproof2}{Part (i)}
    
We observe that, by the definition of  proximal gradient mapping and the triangle inequality, the following holds
\begin{eqnarray}
    && \norm{ G_{t}^{\mu} (x_{t}) } = \mu \norm{ S_{t}^{\mu} (x_t) - x_t } \nonumber \\ 
    &\leq& \mu \norm{ S_{t-1}^{\mu} (x_t) - x_t }
+ \mu \norm{ S_{t}^{\mu} (x_t) - S_{t-1}^{\mu} (x_t)  }
= \norm{ G_{t-1}^{\mu} (x_t)  }
+ \mu \norm{ S_{t}^{\mu} (x_t) - S_{t-1}^{\mu} (x_t)  }   \label{eq:decompose_regret_one_step} . 
\end{eqnarray}

In what follows, we proceed to bound $\norm{S_t^\mu(x_t)-S_{t-1}^\mu(x_t)}$. 
To this end, for fixed $x$, we write 
$
S_t^\mu(x)=\arg \min _{z \in \tilde{\mathcal{L}}(x) } 
\varphi_t(z ; x),
$
where
\begin{eqnarray*}
    \varphi_t(z ; x) &\defeq& 
    \frac{1}{w} \sum_{s=t-w+1}^t f_s\left(\theta_s(x)+J \theta_s(x)(z-x)\right)  \\ 
    && -\left\langle \frac{1}{w} \sum_{s=t-w+1}^t J \theta_s(x)^{\top} \nabla g_s\left(\theta_s(x)\right), z-x\right\rangle+\frac{\mu}{2}\|z-x\|^2  . 
\end{eqnarray*}
We recall that for any $\mu$-strongly convex function $F$ with minimizer $z^*$, the following inequality holds: 
$$
F(z)-F\left(z^*\right) \geq \frac{\mu}{2}\left\|z-z^*\right\|^2 .
$$
Applying this fact twice to $\varphi_{t}(\cdot;x)$ and $\varphi_{t-1}(\cdot;x)$ respectively, we have  
\begin{eqnarray}
   \varphi_t\left( S_{t-1}^\mu (x) ; x \right)-\varphi_t \left(S_t^\mu(x) ; x \right) &\geq& \frac{\mu}{2}\left\|S_{t-1}^\mu(x)-S_t^\mu(x)\right\|^2 ,  \\
   \varphi_{t-1} \left(S_t^\mu (x) ; x \right) - \varphi_{t-1}\left(S_{t-1}^\mu (x) ; x \right) &\geq& \frac{\mu}{2}\left\|S_t^\mu(x)-S_{t-1}^\mu(x)\right\|^2 . 
\end{eqnarray}
For brevity, we suppress the argument $x$ whenever the context is clear. 
Summing the above two inequalities yields that 
\begin{equation} \label{eq:s_s_varphi}
    \mu\left\|S_t^\mu-S_{t-1}^\mu\right\|^2 
\leq \varphi_t\left(S_{t-1}^\mu\right)-\varphi_{t-1}\left(S_{t-1}^\mu\right)+\varphi_{t-1}\left(S_t^\mu\right)-\varphi_t\left(S_t^\mu\right) . 
\end{equation}
We investigate the right-hand-side of \eqref{eq:s_s_varphi} further. For any $z_1, z_2 \in \mathbb{R}^n$, we have 
\begin{eqnarray}
   && \paren{ \varphi_{t}(z_1) - \varphi_{t-1}(z_1) } - \paren{ \varphi_{t}(z_2) - \varphi_{t-1}(z_2) } \nonumber \\ 
   &=& \frac{1}{w} \Big [ f_{t} \paren{ \theta_{t} (x) + J \theta_{t} (x) (z_1 - x)  } - f_{t-w} \paren{ \theta_{t-w} (x) + J \theta_{t-w} (x) (z_1 - x)  }   \nonumber \\ 
   && - f_{t} \paren{ \theta_{t} (x) + J \theta_{t} (x) (z_2 - x)  } 
   +  f_{t-w} \paren{ \theta_{t-w} (x) + J \theta_{t-w} (x) (z_2 - x)  }   \Big ]    \nonumber \\ 
   & & - \inner{ v_{t} - v_{t-1} , z_1-z_2 } \nonumber \\ 
   &\leq& 2  \frac{1}{w} L_{f} L_{\theta} \norm{z_1-z_2} 
   + 2 \frac{1}{w} L_g L_{\theta}  \norm{z_1-z_2} . 
\end{eqnarray}
The inequality follows from Assumption~\ref{assumption:regular_triple_obj} that $f_t$ is $L_f$-Lipschitz, $g_t$ is $L_g$-Lipschitz and that $ \norm{ J \theta_t } \leq L_{\theta}$. 
Now, taking $z_1=S_{t-1}^{\mu} (x_t)$ and $z_2=S_{t}^{\mu} (x_t)$, we have 
\begin{equation}
    \mu \norm{ S_t^\mu (x_t)-S_{t-1}^\mu (x_t) }^2 \leq 2 \frac{1}{w} (L_f+L_g) L_{\theta} \norm{ S_t^\mu (x_t) -S_{t-1}^\mu (x_t) } . 
\end{equation}
Hence we have  
$
    \mu \norm{ S_t^\mu (x_t)-S_{t-1}^\mu (x_t) } \leq 2 \frac{1}{w} (L_f+L_g) L_{\theta}   . 
$
Continuing from \eqref{eq:decompose_regret_one_step}, we have 
\begin{eqnarray}
        \norm{ G_{t}^{\mu} (x_{t}) }^2 &\leq& \paren{ \norm{ G_{t-1}^{\mu} (x_t)  } 
    + 2 \frac{1}{w} (L_f+L_g) L_{\theta}  }^2 \nonumber \\
    &\leq&  \paren{ \frac{\delta}{w}  + 2 \frac{1}{w} (L_f+L_g) L_{\theta}   }^2  .  
\end{eqnarray}

\end{myproof2}



\begin{myproof2}{Part (ii)}
    
In what follows, we denote $z=z_{t+1}^{k}$ and $z'=z_{t+1}^{k+1}$ for brevity.
First, we make the following observations: 
\begin{itemize}
    \item For any $t$ (and $k$), we have  
    \begin{equation}
        f_t( \theta_t ( z' ) ) 
    - f_t( \theta_t (z) + J \theta_t (z) \paren{z'-z} )  \nonumber  \\ 
    \leq L_f \norm{ \theta_t ( z' ) - \paren{ \theta_t (z) + J \theta_t (z) \paren{z'-z} } } \\ 
    \leq \frac{L_f \beta_{\theta} }{2} \norm{z'-z}^2  . 
    \end{equation}
    The first inequality is due to the Lipschitz continuity of $f_t$ and the second inequality is due to the smoothness of $\theta_t$. Hence, we have 
    \begin{eqnarray}
        & & \frac{1}{w} \sum_{s=t-w+1}^{t} \brackets{ f_s ( \theta_s (z') ) - f_s ( \theta_s (z) + J \theta_s (z) \paren{z'-z} ) } \nonumber \\ 
        &\leq& \frac{1}{w} \sum_{s=t-w+1}^{t}  L_f \norm{  \theta_s (z') - \paren{ \theta_s (z) + J \theta_s (z) \paren{z'-z}  }  } \nonumber \\ 
        &\leq& \frac{L_f \beta_{\theta} }{2} \norm{z'-z}^2  . \label{eq:f_step}
    \end{eqnarray}

    \item 
    Regarding the function $g_t$, we have 
    \begin{eqnarray}
        && g_t( \theta_t (z') )  \nonumber \\ 
        &\geq& g_t(\theta_t (z)) 
        + \inner{ \nabla g_t (\theta_t(z))  ,  \theta_t (z') - \theta_t (z) } \label{eq:g_convexity}  \\ 
        &=& g_t(\theta_t (z)) + 
        \nabla g_t (\theta_t(z))^{\top} 
        \paren{ \theta_t(z') - \theta_t(z) - J \theta_t(z) \paren{z'-z  } }
        + \nabla g_t (\theta_t(z))^{\top} J \theta_t(z) \paren{z'-z  } \label{eq:g_add_subtract}  \\ 
        &\geq& g_t(\theta_t (z)) 
        + \nabla g_t (\theta_t(z))^{\top} J \theta_t(z) \paren{z'-z  } 
        - \norm{ \nabla g_t (\theta_t(z)) }  \norm{  \theta_t(z') - \theta_t(z) - J \theta_t(z) \paren{ z'-z  } }  \label{eq:g_reverse_cauchy} \\ 
        &\geq& g_t(\theta_t (z)) 
        + \nabla g_t (\theta_t(z))^{\top} 
        J \theta_t(z) \paren{ z'-z  } 
        -  \frac{L_g \beta_{\theta}}{2} \norm{ z-z'}^2 .  \label{eq:g_smooth_theta}
    \end{eqnarray}
    Inequality~\eqref{eq:g_convexity} is due to the convexity of $g_t$.
    Equation~\eqref{eq:g_add_subtract} follows by adding and subtracting the same term. 
    Inequality~\eqref{eq:g_reverse_cauchy} is due to the Cauchy-Schwarz inequality. 
    Finally, \eqref{eq:g_smooth_theta} is due to the smoothness of $g_t$.
    Therefore, we have 
    \begin{equation} \label{eq:g_step}
        \frac{1}{w} \sum_{s=t-w+1}^{t} g_s \paren{ \theta_s \paren{ z' } }  \geq  
    \frac{1}{w} \sum_{s=t-w+1}^{t} \brackets{ 
        g_s \paren{ \theta_s \paren{ z } }  
    + \nabla g_s (\theta_s(z))^{\top} 
    J \theta_s(z) \paren{ z'-z  }  }  
    -  \frac{L_g \beta_{\theta}}{2} \norm{ z-z' }^2  . 
    \end{equation}

    \item 
    Let $v_g = \frac{1}{w} \sum_{s=t-w+1}^{t} J \theta_s (z)^{\top}  \nabla g_s ( \theta_s (z) ) $. 
    By design of the algorithm, we have 
    $$
    z' \leftarrow \argmin_{z' \in \tilde{\mathcal{L}}(z) }    \frac{1}{w} \sum_{s=t-w+1}^{t} f_s \paren{ \theta_s (z) +  J \theta_s (z) (z'-z) }  
    - \inner{v_{g}^{} , z'-z}
    +  \frac{\mu}{2  } \norm{ z'-z }^2 ,
    $$
    which implies that (since $z$ is a feasible point as well)
    \begin{equation} \label{eq:optimality_step}
        \frac{1}{w} \sum_{s=t-w+1}^{t} f_s \paren{ \theta_s (z) +  J \theta_s (z) (z'-z) }  
        - \inner{ v_{g}  , z'-z}
        +  \frac{\mu}{2  } \norm{ z'-z }^2
        \leq \frac{1}{w} \sum_{s=t-w+1}^{t} f_s \paren{ \theta_s (z) }    .   
    \end{equation}

\end{itemize}

Now, we are ready to put together all the pieces to quantify the descent of $\Phi_t$ from $z$ to $z'$. 
\begin{eqnarray}
    \Phi_t( z' ) 
    &=& \frac{1}{w} \sum_{s=t-w+1}^{t} \brackets{ f_s( \theta_s ( z') ) - g_s( \theta_s ( z' ) )  }  \nonumber   \\ 
    &\overset{\text{(a)}}{\leq}&   \frac{1}{w} \sum_{s=t-w+1}^{t} \brackets{ f_s ( \theta_s (z) + J \theta_s (z) (z'-z) ) } 
    + \frac{L_f \beta_{\theta} }{2} \norm{ z'-z }^2  
    - \frac{1}{w} \sum_{s=t-w+1}^{t}  g_s( \theta_s ( z' ) ) \nonumber   \\ 
    &\overset{\text{(b)}}{\leq}&  \frac{1}{w} \sum_{s=t-w+1}^{t} f_s \paren{ \theta_s (z) }  
    + \inner{v_g , z'-z}
    - \frac{\mu}{2  } \norm{ z'-z }^2
    + \frac{L_f \beta_{\theta} }{2} \norm{ z'-z }^2  
    - \frac{1}{w} \sum_{s=t-w+1}^{t}  g_s( \theta_s ( z' ) ) \nonumber  \\ 
    &\overset{\text{(c)}}{\leq}&   \frac{1}{w} \sum_{s=t-w+1}^{t} f_s \paren{ \theta_s (z) }  
    + \inner{v_g , z'-z}
    - \frac{\mu}{2  } \norm{ z'-z }^2
    + \frac{L_f \beta_{\theta} }{2} \norm{ z'-z }^2  \nonumber \\
    && - \frac{1}{w} \sum_{s=t-w+1}^{t} \brackets{ g_s \paren{ \theta_s \paren{ z } }  
    + \nabla g_s\left(\theta_s(z)\right)^{\top} 
    J \theta_s(z) (z'-z) 
    }  
    +  \frac{L_g \beta_{\theta}}{2} \norm{ z'-z }^2   \nonumber  \\ 
    &\overset{\text{(d)}}{=}&  \frac{1}{w} \sum_{s=t-w+1}^{t} \brackets{ f_s \paren{ \theta_s (z) } - g_s \paren{ \theta_s \paren{ z } } } 
    + \frac{1}{2} \paren{ (L_f+L_g) \beta_{\theta} - \mu } \norm{z'-z}^2 \nonumber  \\ 
    &=& \Phi_t( z )
    + \frac{1}{2} \paren{ (L_f+L_g) \beta_{\theta} - \mu } \norm{z'-z}^2  . 
\end{eqnarray}
Inequality~(a) follows from Inequality~\eqref{eq:f_step}. 
Inequality~(b) follows from Inequality~\eqref{eq:optimality_step}.
Inequality~(c) follows from Inequality~\eqref{eq:g_step}. 
Equation~(d) is due to simple algebraic manipulation.

By design of the algorithm, the while loop 
$
z' \leftarrow S_{t}^{\mu} (z) 
$
will be executed whenever $ \mu \norm{ S_{t}^{\mu} (z) - z} > \frac{\delta}{w}$, namely, $\norm{z'-z} > \frac{1}{\mu} \frac{\delta}{w}$.
$$
\Phi_t( z' )  - \Phi_t( z ) 
\leq - \frac{1}{2} \paren{ \mu - (L_f+L_g) \beta_{\theta}  } \norm{ z'-z }^2 
< - \frac{1}{2} \frac{\mu - (L_f+L_g) \beta_{\theta}}{\mu^2}  \paren{  \frac{\delta}{w} }^2 . 
$$
Let $\tau_t$ be the number of iterations of the while loop executed  at time $t$. We have 
\begin{equation} \label{eq:descent_bound_phi_x}
    \Phi_t( x_{t+1} )  - \Phi_t( x_{t} ) 
< - \frac{1}{2} \frac{\mu - (L_f+L_g) \beta_{\theta}}{\mu^2}  \paren{  \frac{\delta}{w} }^2 \tau_{t} . 
\end{equation}

To proceed, we write the telescoping sum (understanding that $\Phi_0(\cdot)=0$):
\begin{eqnarray}
   \Phi_T(x_{T+1}) - \Phi_0 (x_1) 
   &=& \sum_{t=1}^{T} \paren{ \Phi_t(x_{t+1}) - \Phi_{t-1}(x_{t}) }  \nonumber  \\ 
   &=& \sum_{t=1}^{T} \paren{ \Phi_t(x_{t+1}) - \Phi_{t}(x_{t}) } + \sum_{t=1}^{T} \paren{ \Phi_{t}(x_{t}) - \Phi_{t-1}(x_{t}) }  \nonumber  . 
\end{eqnarray}
The first sum is controlled by \eqref{eq:descent_bound_phi_x}: 
\begin{equation}
    \sum_{t=1}^{T} \paren{ \Phi_t(x_{t+1}) - \Phi_{t}(x_{t}) } 
    \leq - \frac{1}{2} \frac{\mu - (L_f+L_g) \beta_{\theta}}{\mu^2}  \paren{  \frac{\delta}{w} }^2 \sum_{t=1}^{T}\tau_{t} . 
\end{equation}
For the second sum, we note that 
\begin{equation}
    \Phi_t(x)-\Phi_{t-1}(x)=\frac{1}{w}\left(\ell_t(x)-\ell_{t-w}(x)\right)
    \leq \frac{2M}{w} . 
\end{equation}

Hence, we have
\begin{equation}
    \sum_{t=1}^{T} \tau_{t} 
    \leq  \frac{   M + 2 \frac{1}{w} M T }{ \frac{1}{2} \frac{\mu - (L_f+L_g) \beta_{\theta}}{\mu^2}  \paren{  \frac{\delta}{w} }^2 }
    = \frac{2 M \mu^2 }{\mu-\left(L_f+L_g\right) \beta_\theta } \frac{w^2}{\delta^2} \left(1+\frac{2 T}{w}\right)   . 
\end{equation}

If $(L_f+L_g)\beta_{\theta} \neq 0$, setting $\mu = 2 (L_f+L_g)\beta_{\theta}$ yields that
$$
\sum_{t=1}^{T} \tau_{t}  \lesssim M \left(L_f+L_g\right) \beta_\theta \left(\frac{w}{\delta}\right)^2 \paren{1 + \frac{T}{w}}  . 
$$
\end{myproof2}

\begin{remark}[Choice of $\mu$]
    The condition $\mu > (L_f+L_g)\beta_\theta$ plays the same role in our setting as the step-size restriction in time-smoothed online gradient descent \citep{hazan2017efficient}. In both cases, the parameter must be chosen large enough to dominate the curvature term that appears in the descent analysis. 
    If the quantity $(L_f+L_g)\beta_\theta$ is not known a priori, one may instead choose $\mu$ adaptively by a standard backtracking strategy, increasing it until the proximal-linear surrogate yields sufficient descent. Such adaptive curvature selection is classical in line-search and proximal-linear methods; see, for example, \cite{armijo1966minimization,beck2009fast}.
\end{remark}



\subsection{Implications for offline non-convex non-smooth optimization}
\label{subsec:offline_implications}

We now discuss how the online framework recovers the offline composite DC optimization problem as a special case. Consider the fixed objective $\ell(x) = (f-g) \circ \theta (x) $ which satisfies Assumption~\ref{assumption:regular_triple_obj}. We suppose that the same loss is revealed at every round, i.e., $\ell_t = \ell$ for all $t=1,\cdots,T$.

Under the zero-padding convention for $s \leq 0$, the time-smoothed objective simply becomes $\Phi_t(x)=\frac{1}{w} \sum_{s=t-w+1}^t \ell_s(x)=\ell(x)$ for every $t \geq w$. Hence, the online proximal linear mapping reduces to 
$$
    S^{\mu}_{t}(x) = \arg \min_{z \in \tilde{\mathcal{L}} (x)  } ~ \Big\{ f(\theta(x)+J \theta(x)(z-x))-\left\langle J \theta(x)^{\top} \nabla g(\theta(x)), z-x\right\rangle+\frac{\mu}{2}\|z-x\|^2 \Big \}
$$
and 
$
    G^{\mu}_{t}(x) = \mu\left(S_{t}^{\mu}(x)-x\right) := G^{\mu}_{}(x)  . 
$

For $1 \leq t \leq t' \leq T$, we denote by $\mathcal{U}[t,t']$ the uniform distribution over the set $\{t, t+1, \cdots, t'\}$.
Let $\hat{t} \sim \mathcal{U}[w, T]$ be drawn uniformly from the rounds $w, \ldots, T$.
Thus, we have 
\begin{eqnarray}
    \mathbb{E}_{\hat{t} \sim \mathcal{U}[w, T]} \left\|G^\mu \left(x_{\hat{t}}\right)\right\|^2 
    &=& \frac{1}{T-w+1} \sum_{t=w}^T\left\|G^\mu_t \left(x_t\right)\right\|^2 \nonumber 
    \leq \frac{1}{T-w+1} \sum_{t=1}^T\left\|G_t^\mu\left(x_t\right)\right\|^2 = \frac{1}{T-w+1} \regret_T \ . 
\end{eqnarray}
Applying the regret bound in Theorem~\ref{thm:main}, we have
\begin{equation}
    \regret_T \lesssim \frac{T}{T-w+1} \paren{ \frac{\delta + (L_f+L_g) L_{\theta}}{w}}^2 \nonumber . 
\end{equation}
Given a desired accuracy $\varepsilon>0$, choosing  
$$
\delta:=\left(L_f+L_g\right) \beta_{\theta}, \quad 
w:=2\left(L_f+L_g\right) \beta_{\theta} \sqrt{\frac{2}{\varepsilon}}, \quad 
T:=2 w 
$$
yields that $\mathbb{E}_{\hat{t} \sim \mathcal{U}[w, T]} \left\|G^\mu \left(x_{\hat{t}}\right)\right\|^2 \leq \varepsilon$.

\subsection{Computational concerns}
\label{subsec:computational_concerns}





At each time step, we require solving convex optimization problem \eqref{line:trial_d_opt_constrained} multiple times until the stopping criterion is met.
The absence of strict feasibility is well known to cause both theoretical and computational challenges in constrained optimization. 
In what follows, we discuss the condition that guarantees the existence of an interior point for \eqref{line:trial_d_opt_constrained}.

To this end, for a given point $\bar{x}$, we define 
\begin{equation} \label{eq:def_H_feasible_directions}
    H_{\bar{x}} (d) \defeq  \phi \paren{ \zeta(\bar{x}) + J \zeta(\bar{x}) d } 
    - \psi(\zeta(\bar{x}) ) 
    -  \nabla \psi \paren{\zeta(\bar{x})}^{\top} J \zeta ( \bar{x} ) d 
    + \frac{1}{2}  \paren{ L_{\phi} \beta_{\zeta} + L_{\psi} \beta_{\zeta}  } \norm{d}^2 \ . 
\end{equation}
One can easily see that 
\begin{eqnarray*}
    \tilde{\mathcal{L}}  (\bar{x}) 
    &=&  \Bigg\{ \bar{x} + d \in \mathbb{R}^n\,\, \Bigg|\,\,
    \phi \paren{ \zeta(\bar{x}) + J \zeta(\bar{x}) d } 
    + \frac{1}{2} L_{\phi} \beta_{\zeta}  \norm{d}^2 
    \leq \psi(\zeta(\bar{x}) ) +  \nabla \psi \paren{\zeta(\bar{x})}^{\top} J \zeta ( \bar{x} ) d 
    - \frac{1}{2} L_{\psi} \beta_{\zeta}  \norm{d}^2 
    \Bigg\}  \\ 
    &=& \bar{x} + \braces{ d \in \mathbb{R}^n \,\, \Big|\,\, H_{\bar{x}} (d) \leq 0 } .
\end{eqnarray*}
The set $\tilde{\mathcal{L}}  (\bar{x}) $ has nonempty interior if and only if there exists some $d$ such that $H_{\bar{x}} (d) < 0$.
We recall a standard fact in convex analysis, stated in the lemma below. 
\begin{lemma} \label{lemma:interior_iff_condition}
    Let $H: \mathbb{R}^n \rightarrow (-\infty, \infty]$ be a proper convex function. Assume that $H(0) < \infty$. 
    Then we conclude that $0 \notin \partial H (0)$ if and only if there exists $  d \in \mathbb{R}^n$ such that $H(d)<H(0)$.
\end{lemma}
By Lemma~\ref{lemma:interior_iff_condition}, that $\tilde{\mathcal{L}}  (\bar{x}) $ has nonempty interior is equivalent to $0 \notin \partial H_{\bar{x}} (0)$, namely, 
$
    0 \notin J \zeta(\bar{x})^{\top} \paren{ \partial \phi (\zeta(\bar{x})) - \nabla \psi ( \zeta(\bar{x}) ) } .
$
Furthermore, the condition $0 \notin  \partial H_{\bar{x}} (0)$ is equivalent to $\operatorname{dist}\paren{0,  \partial H_{\bar{x}} (0)} > 0$.
Let $\mathcal{Z}$ be a compact set containing all possible iterates $\braces{z_t^k}_{t}^{k}$ of the algorithm. 
Ideally, we want to have 
\begin{equation} \label{eq:want_condition_slater}
    \inf_{ \bar{x} \in \mathcal{Z} \cup \partial \mathcal{X} } \operatorname{dist}\paren{0,  \partial H_{\bar{x}} (0)} > 0 , 
\end{equation}
where the infimum is taken over $\mathcal{Z}$ and the boundary of the domain $\mathcal{X}$, since if $c(\bar{x}) < 0$, then $\bar{x}$ is in the interior of $\mathcal{X}$, and hence $\tilde{\mathcal{L}}  (\bar{x})$ has nonempty interior trivially.

We provide a sufficient condition in the following proposition. 

\begin{proposition} \label{prop:uniform_slater_sufficient_condition}
    A sufficient condition for $\inf_{ \bar{x} \in \mathcal{Z} \cup \partial \mathcal{X} } \operatorname{dist}\paren{0,  \partial H_{\bar{x}} (0)} > 0$ to hold is that 
    \begin{equation}
        \inf_{ \bar{x} \in \mathcal{Z} \cup \partial \mathcal{X} } ~ \braces{  \sigma_{\min}^{+}  \paren{ J \zeta(\bar{x}) }  
        \cdot \operatorname{dist} \paren{ \partial \phi (\zeta(\bar{x})) - \nabla \psi ( \zeta(\bar{x}) ) , \operatorname{null}( J \zeta(\bar{x})^{\top} ) } } > 0 . 
    \end{equation}
\end{proposition}

\subsection{ Error bound in terms of distance to stationarity }
\label{subsec:error_bound}

At time $t$, another natural candidate for measuring the stationarity is $\operatorname{dist}\left(0, \widehat{\partial} \Phi_t(\cdot)+\mathcal{T}_{\mathcal{X}}(\cdot)^*\right)$.
This quantity directly captures how far a point is from satisfying a first-order stationarity condition for the composite constrained problem. 
However, this quantity is practically challenging to compute. 
In contrast, the residual mapping $\norm{G_{t}^{\mu} (\cdot) }$ is more computationally tractable. The goal of this section is to relate these two quantities.
To be specific, given a feasible point $x \in \mathcal{X}$ and $ x^+ = S_{t}^{\mu} (x)$, we provide sufficient conditions under which $\norm{G_{t}^{\mu} (x) }$ can be used to upper bound $\operatorname{dist}\left(0, \widehat{\partial} \Phi_t(x^+)+\mathcal{T}_{\mathcal{X}}(x^+)^*\right)$. 
Consequently, when the algorithm terminates with a small proximal linear residual, we are able to affirmatively assert that the iterate at termination is indeed an approximate stationary point.

\begin{proposition} \label{prop:finite_error_bound}
    Let $x \in \mathcal{X}$ be a feasible point and $x^+ = S_{t}^{\mu} (x)$ be the output of the proximal linear mapping at time $t$ with input $x$. 
    We assume that 
    \begin{itemize}
        \item Suppose Assumption~\ref{assumption:regular_triple_obj} holds.
        \item At time $t$, for every $s \in \{t-w+1, \cdots, t\}$, we assume that there exists some Lagrange multiplier $\lambda \geq 0$ and 
        \begin{equation}
            p_s \in \partial f_s \left( \theta_s (x) + J \theta_s (x) (x^+-x) \right) 
            \quad \text{and} \quad
            q \in \partial \phi \left( \zeta(x) + J \zeta(x) (x^+-x) \right) 
        \end{equation}
        such that the KKT conditions hold with $(\lambda, \braces{p_s}_{s=t-w+1}^t, q)$ for the proximal linear subproblem.  
        Moreover, there exists 
    \begin{equation}
        p_s^+ \in \partial f_s \left( \theta_s (x^+) \right)
         \quad \text{and} \quad
        q^+ \in \partial \phi \left( \zeta(x^+) \right) ,  \nonumber
    \end{equation}
    such that there exists some $\kappa_f , \kappa_{\phi} > 0$ with 
    \begin{equation} \label{eq:assumption_p_s}
        \norm{p_s - p_s^+} \leq \kappa_f \norm{ \theta_s\left(x^{+}\right)-\left(\theta_s(x)+J \theta_s(x) (x^+-x)\right)  }
    \end{equation}
    and 
    \begin{equation} \label{eq:assumption_q}
        \norm{q - q^+} \leq \kappa_{\phi} \norm{ \zeta(x^+) - \left( \zeta(x) + J \zeta(x) (x^+-x) \right) } . 
    \end{equation}

        \item There exists $D > 0$, $\eta > 0$ and a vector $u_x \in \mathbb{R}^n$ such that $\norm{u_x} \leq D$ and $H_x(u_x) \leq -\eta$. 
        \item Assume that  $\lambda > 0$ implies that $x^+ \in \partial \mathcal{X}$, and the Slater condition is satisfied at $x^+$, namely $\exists ~ d^0 \in \mathbb{R}^n \text { with } \phi^{\prime}\left(\zeta\left(x^{+}\right) ; J \zeta\left(x^{+}\right) d^0\right)<\nabla \psi\left(\zeta\left(x^{+}\right)\right)^{\top} J \zeta\left(x^{+}\right) d^0$. 
    \end{itemize}
 
    Then we conclude that 
    \begin{small}
    \begin{eqnarray}
        && \operatorname{dist}\left(0, \widehat{\partial} \Phi_t(x^+)+\mathcal{T}_{\mathcal{X}}(x^+)^*\right) \nonumber \\ 
        &\leq&  \brackets{ \frac{1}{\mu}\left(\beta_\theta\left(L_f+L_g\right)+L_\theta^2 \beta_g+\mu\right) 
        +  +\frac{1}{\mu \eta}\left(2 \beta_\zeta\left(L_\phi+L_\psi\right)+L_\zeta^2 \beta_\psi\right)\left(L_\theta\left(L_f+L_g\right)+\left\|G_t^\mu(x)\right\|\right)\left(D+\frac{\left\|G_t^\mu(x)\right\|}{\mu}\right)
        }  \left\|G_t^\mu(x)\right\| \nonumber \\
        && \brackets{ \frac{1}{2 \mu^2} L_\theta \kappa_f \beta_\theta + +\frac{1}{2 \mu^2 \eta} L_\zeta \kappa_\phi \beta_\zeta\left(L_\theta\left(L_f+L_g\right)+\left\|G_t^\mu(x)\right\|\right)\left(D+\frac{\left\|G_t^\mu(x)\right\|}{\mu}\right) } \left\|G_t^\mu(x)\right\| ^2 . 
    \end{eqnarray}   
    \end{small}
    
    In particular, there exists some constant $\rho > 0$ such that when $\norm{x^+-x} \leq \rho$, we have 
    \begin{equation}
        \operatorname{dist}\left(0, \widehat{\partial} \Phi_t(x^+)+\mathcal{T}_{\mathcal{X}}(x^+)^*\right) 
        \lesssim \left\|G_t^\mu(x)\right\| . 
    \end{equation}

\end{proposition}

The assumptions of Proposition~\ref{prop:finite_error_bound} serve to control the gap between the linearized subproblem \eqref{line:trial_d_opt_constrained} and the first-order stationarity measure evaluated at the updated point $x^+$. 
The conditions \eqref{eq:assumption_p_s} and \eqref{eq:assumption_q} ensure that the subgradient information at the subgradient information at $x$ can be used to inform the subgradient information at the true updated point $x^+$.
The existence of $u_x$ with $H_x(u_x) \leq -\eta$ is used to control the size of the Lagrange multiplier associated with the problem \eqref{line:trial_d_opt_constrained}. 
The condition that $\lambda > 0$ implies $x^+ \in \partial \mathcal{X}$ indicates that when the surrogate constraint is active at $x$, the updated point $x^+$ will be on the boundary of the true domain $\mathcal{X}$.

\section{Concluding remarks}

In this paper, we study online optimization for a broad class of structured non-convex non-smooth problems with both CDCS objectives and CDCS constraints. We propose a proximal-linear local-regret framework and show that the associated residual provides a proper and computationally tractable measure of stationarity for the original problem. Our analysis relies on a tangent-cone characterization for feasible regions described by CDCS constraints, which is of independent interest, and leads to local-regret, iteration-complexity, and error-bound guarantees. 


%
%
%

\begin{APPENDICES}
    \section{Omitted Proofs} 

    \subsection{Proof of Lemma~\ref{lemma:max_of_dc_is_dc}}

\begin{myproof} 
    We stack the two mappings $\theta_1$ and $\theta_2$ to form a new mapping $\theta(x) = (\theta_1(x), \theta_2(x))$, which is differentiable. 
    On the product space, we define 
    \begin{equation}
        f\left(y_1, y_2\right)=\max \left\{f_1\left(y_1\right)+g_2\left(y_2\right), f_2\left(y_2\right)+g_1\left(y_1\right)\right\} ,  \nonumber
    \end{equation}
    which is convex in $(y_1, y_2)$. We also define 
    \begin{equation}
    g\left(y_1, y_2\right)=g_1\left(y_1\right)+g_2\left(y_2\right) , \nonumber
    \end{equation}
    which is convex and differentiable in $(y_1, y_2)$.

    We note that 
    \begin{eqnarray}
        & & \max \braces{ \phi_1(x) , \phi_2(x) } \nonumber \\ 
        &=& \max \braces{ f_1(\theta_1(x)) - g_1(\theta_1(x)) , f_2(\theta_2(x)) - g_2(\theta_2(x)) } \nonumber \\
        &=& \max \braces{ 
            f_1\left(\theta_1(x)\right)+g_2\left(\theta_2(x)\right), 
            f_2\left(\theta_2(x)\right)+g_1\left(\theta_1(x)\right) 
            } 
        - \left(g_1\left(\theta_1(x)\right)+g_2\left(\theta_2(x)\right)\right) \label{eq:max_constraints_using_basic_algebra} \\
        &=& f(\theta(x))-g(\theta(x)) . 
    \end{eqnarray}
    The last equality follows from the definitions of $f, g, \theta$. 
    Equation \eqref{eq:max_constraints_using_basic_algebra} holds by noticing the fact that for any real numbers $a, b, c, d$, we have 
    \begin{equation} 
        \max \braces{ a-b , c-d } = \max \braces{ a + d , b + c } - (b+d) . \nonumber
    \end{equation}

\end{myproof}

    \subsection{Proof of Lemma~\ref{lemma:zero_in_subdiff_implies_bouligand_stationarity}}
\begin{myproof}
    The fact $0 \in \widehat{\partial} \phi (\bar{x}) + \mathcal{T}_{X} ( \bar{x} )^*$ implies that 
    \begin{equation}
        0=g+v \text { with } g \in \widehat{\partial} \phi(\bar{x}) \text { and } v \in \mathcal{T}_X(\bar{x})^*
    \end{equation}
    where $\inner{v,d} \leq 0 \text { for all } d \in \mathcal{T}_X(\bar{x})$.

    By the supporting-hyperplane inequality for the regular subdifferential, we have 
    \begin{eqnarray}
        \phi'( \bar{x}; d ) &\geq& \sup_{ s \in \widehat{\partial} \phi(\bar{x}) } ~ \inner{ s , d } 
        \geq \inner{g, d}  
        = - \inner{v, d}  
        \geq 0 , 
    \end{eqnarray}
    for all $d \in \mathcal{T}_{X} ( \bar{x} )$. 
\end{myproof}

\subsection{Proof of Lemma~\ref{lemma:interior_iff_condition}}
\begin{myproof}
    By definition of the convex subdifferential,
$$
v \in \partial H(0) \quad \Longleftrightarrow \quad H(d) \geq H(0)+\langle v, d-0\rangle \text { for all } d \in \mathbb{R}^n .
$$
Now set $v=0$. Then
$$
0 \in \partial H(0) \quad \Longleftrightarrow \quad H(d) \geq H(0) \text { for all } d \in \mathbb{R}^n .
$$
Taking negation gives the desired claim. 
\end{myproof}


\subsection{Proof of Proposition~\ref{prop:uniform_slater_sufficient_condition}}

To proceed, we first state some standard results in linear algebra. 
We include the proofs for completeness.
\begin{lemma} \label{lemma:lower_bound_Ax_by_svd}
   For any matrix $A \in \mathbb{R}^{m \times n}$ and any vector $w \in \mathbb{R}^m$, we have 
    \begin{equation}
        \left\|A^{\top} w \right\| \geq \sigma_{\min }^{+}(A)\left\|\mathcal{P}_{\operatorname{range}(A)} w \right\| . 
    \end{equation}
\end{lemma}

\begin{myproof}
    If $\operatorname{rank}(A) = 0$, then $A^{\top} w = 0$ for all $w$, and the inequality holds trivially.
    Now suppose $\operatorname{rank}(A) = r > 0$.
    Let the singular value decomposition of $A$ be $A = U \Sigma V^{\top}$, where $U \in \mathbb{R}^{m \times m}$ and $V \in \mathbb{R}^{n \times n}$ are orthonormal matrices, and $\Sigma \in \mathbb{R}^{m \times n}$ is a rectangular diagonal matrix with singular values $\sigma_1 \geq \sigma_2 \geq \cdots \geq \sigma_r > 0$ on the diagonal and zeros elsewhere.
    Hence, 
    $$
    \norm{A^\top w } = \norm{ V \Sigma^{\top} U^{\top} w } = \norm{ \Sigma^{\top} U^{\top} w } , 
    $$
    where the second equality follows from the orthonormality of $V$.
    Let $y = U^{\top} w$. Then 
    $$
    \norm{ \Sigma^{\top} U^{\top} w }
    = \sqrt{ \sum_{i=1}^{r} \sigma_i^2 y_i^2 } 
    \geq \sigma_r \sqrt{ \sum_{i=1}^{r} y_i^2 } . 
    $$
    We recall that $\operatorname{range}(A) = \operatorname{span} \braces{ u_1, \cdots, u_r }$, and $\{u_i\}$ are orthonormal columns of $U$.
    Hence, the projection of $w$ onto $\operatorname{range}(A)$ is 
    $$
    \mathcal{P}_{\operatorname{range}(A)} w
    = U\left[\begin{array}{ll}
I_r & \\
& 0
\end{array}\right] U^{\top} w . 
    $$
    Therefore, using again $U$ is orthogonal, $\norm{\mathcal{P}_{\operatorname{range}(A)} w} 
    = \norm{U\left[\begin{array}{ll}
    I_r & \\
    & 0
    \end{array}\right] U^{\top} w}
    = \norm{\left[\begin{array}{ll}
    I_r & \\
    & 0
    \end{array}\right] y } = \sqrt{ \sum_{i=1}^{r} y_i^2 }$.

\end{myproof}

\begin{lemma} \label{lemma:decomposition_distance}
    Let $V$ be a finite-dimensional inner product space, and let $S \subseteq V$ be a subspace. Then for every $v \in V$,
$$
\left\|\mathcal{P}_S v\right\|=\operatorname{dist}\left(v, S^{\perp}\right) = \inf _{w \in S^{\perp}}\|v-w\| .
$$
\end{lemma}

\begin{myproof}

    Since $V=S \oplus S^{\perp}$ orthogonally, every $v \in V$ admits the unique decomposition
$
v=\mathcal{P}_S v+\mathcal{P}_{S^{\perp}} v,
$
with $\mathcal{P}_S v \in S$ and $\mathcal{P}_{S^{\perp}} v \in S^{\perp}$.
Now we take any $w \in S^{\perp}$. Then
$$
v-w=\mathcal{P}_S v+\left(\mathcal{P}_{S^{\perp}} v-w\right),
$$
where
$\mathcal{P}_S v \in S$ , $\mathcal{P}_{S^{\perp}} v-w \in S^{\perp}$.
The former two terms are orthogonal, so we have 
$$
\|v-w\|^2=\left\|\mathcal{P}_S v\right\|^2+\left\|\mathcal{P}_{S^{\perp}} v-w\right\|^2 \geq\left\|\mathcal{P}_S v\right\|^2 .
$$
Hence, we conclude 
$$
\|v-w\| \geq\left\|\mathcal{P}_S v\right\| \quad \forall w \in S^{\perp} .
$$
Taking the infimum over $w \in S^{\perp}$ gives
$$
\operatorname{dist}\left(v, S^{\perp}\right) \geq\left\|\mathcal{P}_S v\right\| .
$$

On the other hand, choosing $w=P_{S^{\perp}} v \in S^{\perp}$, we get
$
\|v-w\|=\left\|v-\mathcal{P}_{S^{\perp}} v\right\|=\left\|\mathcal{P}_S v\right\| .
$
Therefore
$
\operatorname{dist}\left(v, S^{\perp}\right) \leq\left\|\mathcal{P}_S v\right\| .
$

Combining both inequalities, we conclude that 
$$
\operatorname{dist}\left(v, S^{\perp}\right)=\left\|\mathcal{P}_S v\right\| .
$$
The proof of the lemma is now complete.
\end{myproof}

We are now ready to prove Proposition~\ref{prop:uniform_slater_sufficient_condition}. 
In view of Lemma~\ref{lemma:lower_bound_Ax_by_svd}, we have 
\begin{eqnarray}
    \operatorname{dist}\paren{0,   \partial H_{\bar{x}} (0) } 
    &=& \inf_{ q \in \partial \phi (\zeta(\bar{x})) - \nabla \psi ( \zeta(\bar{x}) ) } \norm{ J \zeta(\bar{x})^{\top} q } \nonumber \\
    &\geq& \sigma_{\min}^{+}  \paren{ J \zeta(\bar{x}) } \cdot  \inf_{ q \in \partial \phi (\zeta(\bar{x})) - \nabla \psi ( \zeta(\bar{x}) ) }  \norm{ \mathcal{P}_{\operatorname{range}(J \zeta(\bar{x}))} (q) }   \\ 
    &=&  \sigma_{\min}^{+}  \paren{ J \zeta(\bar{x}) } \cdot \operatorname{dist} \paren{ \partial \phi (\zeta(\bar{x})) - \nabla \psi ( \zeta(\bar{x}) ) , \operatorname{null}( J \zeta(\bar{x})^{\top} ) } \label{eq:apply_decomposition_algebra} . 
\end{eqnarray}
Here we explain the reasoning for \eqref{eq:apply_decomposition_algebra}. 
Because
$
\mathbb{R}^m=\operatorname{range}(J \zeta(\bar{x})) \oplus \operatorname{null}\left(J \zeta(\bar{x})^{\top}\right)
$
is an orthogonal decomposition, by Lemma~\ref{lemma:decomposition_distance}, we have
$$
\left\|\mathcal{P}_{\mathrm{range}(J \zeta(\bar{x}))}(q)\right\|=\operatorname{dist}\left(q, \operatorname{null}\left(J \zeta(\bar{x})^{\top}\right)\right) .
$$
Hence, 
$$
\inf _{q \in  \partial \phi (\zeta(\bar{x})) - \nabla \psi ( \zeta(\bar{x}) )  }\left\|\mathcal{P}_{\text {range }(J \zeta(\bar{x}))}(q)\right\|=\operatorname{dist}\left(  \partial \phi (\zeta(\bar{x})) - \nabla \psi ( \zeta(\bar{x}) )  , \operatorname{null}\left(J \zeta(\bar{x})^{\top}\right)\right) .
$$

\subsection{Proof of Proposition~\ref{prop:finite_error_bound}}

\begin{myproof}
    
For brevity, we denote $ d = x^+ - x$. 
We denote $v = \frac{1}{w}   \underset{s=t-w+1}{\overset{t}{\sum}}  J \theta_s (x)^{\top}  \nabla g_s ( \theta_s (x) )$.
By design of the algorithm in \eqref{line:trial_d_opt_constrained} and the definition of the set $H$ in \eqref{eq:def_H_feasible_directions}, we have that $d$ is the unique minimizer of the following convex problem 
\begin{equation}
    \min_{ u : H_x(u) \leq 0 } ~ \braces{ \frac{1}{w} \sum_{s=t-w+1}^{t} f_s \paren{ \theta_s (x) + J \theta_s (x) u }
    - \inner{  v  , u}
    + \frac{\mu}{2} \norm{u}^2 }  . 
\end{equation}
There exists some Lagrange multiplier $\lambda \geq 0$ such that the KKT conditions hold at the chosen subgradients $p_s$ and $q$ 
\begin{equation}  \label{eq:KKT_conditoin_error_bound}
    0 = \frac{1}{w} \sum_{s=t-w+1}^{t}  J \theta_s (x)^{\top} p_s - v 
    + \mu d 
    + \lambda \cdot \paren{ J \zeta (x)^{\top} \paren{ q - \nabla \psi (\zeta(x))}  + \paren{ L_{\phi} \beta_{\zeta} + L_{\psi} \beta_{\zeta} } d }  , 
\end{equation}
together with 
\begin{equation}
    H_x(d) \leq 0 , \quad \lambda \geq 0 , \quad \lambda \cdot H_x(d) = 0 . 
\end{equation}
Since $H_x(d) \leq 0$, we have $x^+ = x + d \in \tilde{\mathcal{L}}(x)$. By Theorem~\ref{thm:tangent_cone_characterization}, we have $x^+ \in \mathcal{X}$. 
We denote
\begin{eqnarray*}
    a_s^+ &=& J \theta_s ( x^+)^{\top} p_s^+ - J \theta_s( x^+ )^{\top} \nabla g_s ( \theta_s(x^+) ) ,\\
    a^+ &=& \frac{1}{w} \sum_{s=t-w+1}^t a_s^+ ,\\ 
    a &=& \frac{1}{w} \sum_{s=t-w+1}^t J \theta_s (x)^{\top} p_s - v  ,\\
    b^+ &=& J \zeta (x^+)^{\top} \paren{ q^+ - \nabla \psi (\zeta(x^+))} ,\\ 
    b &=& J \zeta (x)^{\top} \paren{ q - \nabla \psi (\zeta(x))}  + \paren{ L_{\phi} \beta_{\zeta} + L_{\psi} \beta_{\zeta} } d . 
\end{eqnarray*}
One can verify that $a^+ \in \widehat{\partial} \Phi_t(x^+)$. 
Then clearly, we have 
\begin{equation}
    \operatorname{dist}\left(0, \widehat{\partial} \Phi_t(x^+)+\mathcal{T}_{\mathcal{X}}(x^+)^*\right) 
    \leq \norm{a^+ + \mathcal{P}_{\mathcal{T}_{\mathcal{X}}(x^+)^*} (\lambda b^+ )} .
\end{equation}
The KKT condition in \eqref{eq:KKT_conditoin_error_bound} can be written as $0 = a + \mu d + \lambda b$ and hence 
$
    a^+ + \mathcal{P}_{\mathcal{T}_{\mathcal{X}}(x^+)^*} (\lambda b^+ ) 
    = a^+ - a + \mathcal{P}_{\mathcal{T}_{\mathcal{X}}(x^+)^*} (\lambda b^+ ) - \lambda b - \mu d    . 
$
Taking norms on both sides and applying the triangle inequality yields that
\begin{equation} \label{eq:decomposition_step_for_error_bound}
    \norm{a^+ + \mathcal{P}_{\mathcal{T}_{\mathcal{X}}(x^+)^*} (\lambda b^+ ) }
    \leq \norm{ a^+ - a } + \norm{ \mathcal{P}_{\mathcal{T}_{\mathcal{X}}(x^+)^*} (\lambda b^+ ) - \lambda b } + \mu \norm{d} .
\end{equation}

We proceed to bound the three terms on the right-hand side separately. 

\begin{itemize}
    \item For the first term in \eqref{eq:decomposition_step_for_error_bound}, by triangle inequality, we have 
    \begin{eqnarray}
        && \norm{ a^+ - a } \nonumber \\ 
        &=& \norm{ \frac{1}{w} \sum_{s=t-w+1}^t \brackets{ J \theta_s ( x^+)^{\top} p_s^+ - J \theta_s(x^+)^{\top} \nabla g_s ( \theta_s(x^+) ) } - \frac{1}{w} \sum_{s=t-w+1}^t \brackets{ J \theta_s (x)^{\top} p_s - v } } \nonumber \\ 
        &\leq& \frac{1}{w}  \sum_{s=t-w+1}^t  \norm{  J \theta_s ( x^+)^{\top} p_s^+ - J \theta_s (x)^{\top} p_s  } 
        + \frac{1}{w}  \sum_{s=t-w+1}^t  \norm{   J \theta_s(x^+)^{\top} \nabla g_s ( \theta_s(x^+) ) - J \theta_s(x)^{\top} \nabla g_s\left(\theta_s(x)\right) }   \nonumber . 
    \end{eqnarray}
    For the first term, we have 
    \begin{eqnarray}
        && \norm{  J \theta_s ( x^+)^{\top} p_s^+ - J \theta_s (x)^{\top} p_s  }  \nonumber \\ 
        &=& \norm{ \left(J \theta_s\left(x^{+}\right)^{\top}-J \theta_s(x)^{\top}\right) p_s^{+}+J \theta_s(x)^{\top}\left(p_s^{+}-p_s\right) } \nonumber \\
        &\overset{\text{(a)}}{\leq}& \left\|J \theta_s\left(x^{+}\right)-J \theta_s(x)\right\|\left\|p_s^{+}\right\|+\left\|J \theta_s(x)\right\|\left\|p_s^{+}-p_s\right\| \nonumber \\
        &\overset{\text{(b)}}{\leq}& \beta_\theta\|d\| L_f + L_{\theta} \kappa_f \norm{ \theta_s\left(x^{+}\right)-\left(\theta_s(x)+J \theta_s(x) (x^+-x)\right)  }  \nonumber  \\ 
        &\overset{\text{(c)}}{\leq}& \beta_\theta L_f  \|d\|  + \frac{1}{2} L_{\theta} \kappa_f \beta_{\theta} \|d\|^2  . \nonumber 
    \end{eqnarray}
    Inequality~(a) follows from the triangle inequality and the submultiplicativity of the operator norm. 
    For Inequality~(b), we use three facts. First, since $\theta_s$ is $\beta_{\theta}$-smooth, we have $\left\|J \theta_s\left(x^{+}\right)-J \theta_s(x)\right\| \leq \beta_\theta\|d\|$. Second, since $f_s$ is $L_f$-Lipschitz and $p_s^+ \in \partial f_s \left( \theta_s (x^+) \right)$, every subgradient of $f_s$ has norm at most $L_f$, and hence $\norm{p_s^+} \leq L_f$. Third, we recall that we assume the condition in \eqref{eq:assumption_p_s}. 
    Lastly, Inequality~(c) follows from the assumption that $\theta_s$ is $\beta_{\theta}$-smooth.

    Following the same line of reasoning, we can bound  
    \begin{eqnarray}
        & & \left\|J \theta_s\left(x^{+}\right)^{\top} \nabla g_s\left(\theta_s\left(x^{+}\right)\right)-J \theta_s(x)^{\top} \nabla g_s\left(\theta_s(x)\right)\right\| \nonumber  \\
        &\leq&\left\|J \theta_s\left(x^{+}\right)-J \theta_s(x)\right\|\left\|\nabla g_s\left(\theta_s\left(x^{+}\right)\right)\right\|+\left\|J \theta_s(x)\right\|\left\|\nabla g_s\left(\theta_s\left(x^{+}\right)\right)-\nabla g_s\left(\theta_s(x)\right)\right\| \nonumber  \\
        &\leq& \beta_\theta L_g\|d\|+L_\theta^2 \beta_g\|d\| . \nonumber 
    \end{eqnarray}

    Putting the pieces together, we have
    \begin{equation} \label{eq:a_a_plus_upper_bound}
        \norm{ a^+ - a } \leq 
        \paren{ \beta_\theta\left(L_f+L_g\right) + L_\theta^2 \beta_g } \|d\| 
        + \frac{1}{2} L_{\theta} \kappa_f \beta_{\theta} \|d\|^2  .
    \end{equation}

    \item For the second term in \eqref{eq:decomposition_step_for_error_bound}, we can further decompose it as 
    \begin{equation}
        \norm{ \mathcal{P}_{\mathcal{T}_{\mathcal{X}}(x^+)^*} (\lambda b^+ ) - \lambda b }  
        \leq \norm{ \mathcal{P}_{\mathcal{T}_{\mathcal{X}}(x^+)^*} (\lambda b^+ ) - \lambda b^+ } + \norm{ \lambda b^+ - \lambda b }
    = \operatorname{dist}\paren{ \lambda b^+ , \mathcal{T}_{\mathcal{X}}(x^+)^* } + \lambda \norm{ b^+ - b } . 
    \end{equation}
    The last equality directly follows from the definition. 
    
    To proceed, we first provide an upper bound on $\lambda$. 
        If $\lambda = 0$, then there is nothing to prove. 
    Suppose $\lambda > 0$, then the complementary slackness condition implies that $H_x(d) = 0$.
    Since $b \in \partial H_x(d)$, the convex subgradient inequality implies that 
    $
        H_x (u) \geq H_x(d) + \inner{b, u-d} 
    $
    for every $u$. In particular, taking $u = u_x$ yields that 
    \begin{equation}
        \inner{b, d - u_x} \geq \eta . 
    \end{equation}
    Taking inner product with the KKT condition with $d-u_x$ yields that 
    \begin{equation}
        \lambda \eta \leq \inner{ b, d-u_x} \leq \norm{ \frac{1}{w} \sum_{s=t-w+1}^t J \theta_s(x)^{\top} p_s-v+\mu d  } \norm{ d-u_x } . \nonumber 
    \end{equation}
    Since $f_s$ is $L_f$-Lipschitz, $g_s$ is $L_g$-Lipschitz, and $\norm{J \theta_s} \leq L_{\theta}$, we have $ \norm{ \frac{1}{w} \sum_{s=t-w+1}^t J \theta_s(x)^{\top} p_s-v } \leq L_\theta\left(L_f+L_g\right)   $. By triangle inequality, we also have $\norm{d-u_x} \leq \norm{d}+\norm{u_x} \leq \norm{d}+D$.
    Hence, we conclude that 
    \begin{equation}
        \lambda \leq \frac{\left(L_\theta\left(L_f+L_g\right)+\mu\|d\|\right)\left(\|d\|+D\right)}{\eta} . 
    \end{equation}
    
    Next, we provide an upper bound on $\norm{b^+-b}$. 
    \begin{eqnarray}
       \norm{b^+-b} &=& \norm{ J \zeta (x^+)^{\top} \paren{ q^+ - \nabla \psi (\zeta(x^+))} 
       - J \zeta (x)^{\top} \paren{ q - \nabla \psi (\zeta(x))}  
       - \paren{ L_{\phi} \beta_{\zeta} + L_{\psi} \beta_{\zeta} } d
        } \nonumber \\
       &\leq&  \norm{  J \zeta (x^+)^{\top} q^+ -  J \zeta (x)^{\top} q  }
       + \norm{ J \zeta\left(x^{+}\right)^{\top} \nabla \psi\left(\zeta\left(x^{+}\right)\right)-J \zeta(x)^{\top} \nabla \psi(\zeta(x)) }
       +   \paren{ L_{\phi}  + L_{\psi}  } \beta_{\zeta} \norm{d} \nonumber ,
    \end{eqnarray}
    where the inequality follows from the triangle inequality and the definition of $b$ and $b^+$. 

    We first bound the term involving $q^{+}$and $q$. By adding and subtracting $J \zeta(x)^{\top} q^{+}$, we obtain
    \begin{equation}
        \norm{  J \zeta (x^+)^{\top} q^+ -  J \zeta (x)^{\top} q  }
        \leq \left\|J \zeta\left(x^{+}\right)-J \zeta(x)\right\|\left\|q^{+}\right\|+\|J \zeta(x)\|\left\|q^{+}-q\right\| 
        \nonumber . 
    \end{equation}
    where the last step again uses the triangle inequality together with submultiplicativity of the operator norm. Since $J \zeta$ is $\beta_\zeta$-Lipschitz and $\phi$ is $L_\phi$-Lipschitz, we have
    $$
    \left\|J \zeta\left(x^{+}\right)-J \zeta(x)\right\| \leq \beta_\zeta\|d\|, \quad\left\|q^{+}\right\| \leq L_\phi,
    $$
    and hence
    \begin{equation} 
        \left\|J \zeta\left(x^{+}\right)-J \zeta(x)\right\|\left\|q^{+}\right\| \leq \beta_\zeta L_\phi \|d\| . \nonumber
    \end{equation}
    Moreover, by the assumption in \eqref{eq:assumption_q} and the $\beta_\zeta$-smoothness of $\zeta$, we have
    \begin{equation}
        \left\|q^{+}-q\right\| \leq \kappa_\phi\left\|\zeta\left(x^{+}\right)-\zeta(x)-J \zeta(x) d\right\| \leq \frac{1}{2} \kappa_\phi \beta_\zeta\|d\|^2 \nonumber . 
    \end{equation}

    We next bound the term involving $\nabla \psi$. By adding and subtracting $J \zeta(x)^{\top} \nabla \psi\left(\zeta\left(x^{+}\right)\right)$, we have
    \begin{eqnarray}
& & \left\|J \zeta\left(x^{+}\right)^{\top} \nabla \psi\left(\zeta\left(x^{+}\right)\right)-J \zeta(x)^{\top} \nabla \psi(\zeta(x))\right\| \nonumber \\
&\leq& \left\|J \zeta\left(x^{+}\right)-J \zeta(x)\right\|\left\|\nabla \psi\left(\zeta\left(x^{+}\right)\right)\right\|+\|J \zeta(x)\|\left\|\nabla \psi\left(\zeta\left(x^{+}\right)\right)-\nabla \psi(\zeta(x))\right\| \nonumber \\ 
&\leq& \beta_\zeta L_\psi\|d\|+L_{\zeta}^2 \beta_\psi\|d\| . \nonumber
    \end{eqnarray}
    The first inequality again follows from the triangle inequality and the submultiplicativity of the operator norm. The second inequality follows from the $\beta_\zeta$-smoothness of $\zeta$ and the $L_\psi$-Lipschitzness and $\beta_\psi$-smoothness of $\psi$.

    Putting together, we have 
    \begin{equation}
        \norm{b^+-b} \leq \left(2 \beta_\zeta\left(L_\phi+L_\psi\right)+L_\zeta^2 \beta_\psi\right)\|d\|+\frac{1}{2} L_\zeta \kappa_\phi \beta_\zeta\|d\|^2 . 
    \end{equation}

    \item If $\lambda = 0$, then $\operatorname{dist}\paren{ \lambda b^+ , \mathcal{T}_{\mathcal{X}}(x^+)^* } = 0$.
    Otherwise, if $\lambda > 0$, which implies that $x^+ \in \partial \mathcal{X}$ as we assumed, then by Corollary~\ref{corollary:first_order_tangent_cone}, since $q^+ \in \partial \phi ( \zeta(x^+) )$, we have 
    \begin{equation}
        b^{+}=J \zeta\left(x^{+}\right)^{\top}\left(q^{+}-\nabla \psi\left(\zeta\left(x^{+}\right)\right)\right) \in \mathcal{T}_{\mathcal{X}}\left(x^{+}\right)^* . 
    \end{equation}
    Because $\mathcal{T}_{\mathcal{X}}\left(x^{+}\right)^*$ is a cone and $\lambda \geq 0$, we also have $\lambda b^{+} \in \mathcal{T}_{\mathcal{X}}\left(x^{+}\right)^*$. Hence $\operatorname{dist}\paren{ \lambda b^+ , \mathcal{T}_{\mathcal{X}}(x^+)^* } = 0$.
 
\end{itemize}

\vspace{0.5cm}

Putting everything together yields the desired result.

\end{myproof}

\end{APPENDICES}


\bibliographystyle{plainnat}
\bibliography{ref.bib}

\end{document}

%% file: macros.tex
\usepackage{enumitem}
\setlist{leftmargin=*}

\usepackage{algpseudocode}

\usepackage{natbib}
\usepackage[dvipsnames]{xcolor}

\usepackage[hidelinks]{hyperref}

\usepackage{graphicx}
\usepackage{float}
\usepackage{subfigure}

\usepackage{booktabs}

\newcommand{\regret}{\textsf{Regret}} 

\newcommand{\defeq}{\mathrel{\overset{\mathrm{def}}{=}}}

\newcommand{\norm}[1]{\left \| {#1} \right \|}

\newcommand{\inner}[1]{\left \langle {#1} \right \rangle }
\newcommand{\one}[1]{\BFI \left[ {#1} \right]}
\newcommand{\paren}[1]{ \left( {#1} \right) }

\newcommand{\abs}[1]{ \left | {#1} \right | }
\newcommand{\brackets}[1]{\left[ {#1} \right]}
\newcommand{\braces}[1]{\left\{{#1}\right\}}

\newenvironment{myproof}{ \paragraph{Proof: } } {\hfill$\square$}
\newenvironment{myproof2}[1]{ \paragraph{Proof of #1: } } {\hfill$\square$}

